\documentclass[pdflatex,sn-mathphys-ay]{sn-jnl}

\usepackage{graphicx}%
\usepackage{multirow}%
\usepackage{amsmath,amssymb,amsfonts}%
\usepackage{amsthm}%
\usepackage{mathrsfs}%
\usepackage[title]{appendix}%
\usepackage{xcolor}%
\usepackage{textcomp}%
\usepackage{manyfoot}%
\usepackage{booktabs}%
\usepackage{algorithm}%
\usepackage{algorithmicx}%
\usepackage{algpseudocode}%
\usepackage{listings}%

\theoremstyle{thmstyleone}%
\newtheorem{theorem}{Theorem}
\newtheorem{proposition}[theorem]{Proposition}%

\theoremstyle{thmstyletwo}%

\theoremstyle{thmstylethree}%

\newcommand{\dint}{\displaystyle\int}
\newcommand{\dsum}{\displaystyle\sum}
\newcommand{\tht}{\boldsymbol{\theta}}
\newcommand{\thzero}{\boldsymbol{\theta}_0}
\newcommand{\whth}{\widehat{\boldsymbol{\theta}}}
\newcommand{\R}{\mathbb{R}}
\newcommand{\D}{\mathcal{D}}
\newcommand{\ind}{\boldsymbol{1}}
\newcommand{\dd}{\,d}

\def\doiurl#1{} 

\begin{document}

\title[Minimum cumulative power divergence estimation with applications to copulas]{Minimum cumulative power divergence estimation with applications to copulas}


%
%
%

\author*[1]{\fnm{Ostap} \sur{Okhrin}}
\email{ostap.okhrin@tu-dresden.de}

\author[2]{\fnm{Konstantinos} \sur{Zografos}}
\email{kzograf@uoi.gr}

\affil*[1]{%
  \orgdiv{Chair of Applied Statistics,
    Faculty of Transport and Traffic Sciences},
  \orgname{Technische Universit\"at Dresden},
  \orgaddress{%
    \city{Dresden},
    \country{Germany}}}

\affil[2]{%
  \orgdiv{Department of Mathematics},
  \orgname{University of Ioannina},
  \orgaddress{%
    \city{Ioannina},
    \country{Greece}}}


\abstract{Parametric inference based on distribution functions is useful when a model is naturally specified through its cumulative distribution function or when a tractable density is unavailable. We study minimum cumulative power divergence estimation, obtained by comparing a parametric distribution function with its empirical counterpart. For a fixed positive tuning parameter and a fixed finite integration measure, we establish the existence of a consistent sequence of local solutions to the estimating equations and derive its asymptotic normality. The resulting covariance has a sandwich form whose components involve derivatives of the distribution function and the covariance kernel of the multivariate empirical process. We specialize the construction to copulas using Lebesgue measure on the unit cube, and distinguish estimation with known margins from estimation based on ranks. In the latter case, the covariance incorporates the effect of estimating the marginal distributions. The relationship with the copula divergences of Eguchi and Kato is made explicit: the power loss is closely related, whereas the integration measure determines a different estimation criterion. A simulation design for Gumbel, Clayton, and Frank copulas examines the estimation of Kendall's tau under correct specification, family misspecification, and contamination by independence. The simulations show competitive clean-data accuracy but no general robustness advantage over the empirical-copula-weighted alternatives under independence contamination.}




\maketitle

\section{Introduction}\label{sec:intro}

Likelihood inference requires a density or probability mass function that can be evaluated and optimized. This requirement can be inconvenient for multivariate models specified through a cumulative distribution function (CDF), and it excludes a direct density-based treatment of distributions with singular components. CDF offer a different starting point: the empirical CDF estimates them without density smoothing, and a parametric model can be fitted by minimizing a discrepancy between the empirical and model CDFs. The choice of discrepancy then determines which features of the distribution are emphasized by the estimator.

Minimum distance estimation has a long history, including the work of \citet{wolfowitz1957} and \citet{beran1977}. A particularly influential development is the density power divergence of \citet{basu1998}. Its tuning parameter allows the contribution of different density levels to be varied, and the corresponding estimators provide an established framework for studying the balance between efficiency and robustness; see also \citet{basu2011}. Transferring such constructions from densities to distribution or survival functions leads to cumulative information measures. The cumulative residual entropy of \citet{rao2004} and the survival entropies of \citet{zografosnadarajah2005} are early contributions to this literature. \citet{zografos2023} gives a systematic treatment of cumulative and survival counterparts of classical information measures, including a cumulative counterpart of the density power divergence.

The present paper studies the estimation procedure associated with that cumulative power divergence. The criterion compares the parametric CDF directly with the empirical CDF. It therefore requires neither a density estimate nor the existence of a model density. We retain the power parameter $a>0$ and derive the estimating equations, their local consistency and their asymptotic normal distribution. The proofs use the familiar Taylor expansion approach to likelihood and minimum-divergence estimation: the first derivative of the criterion vanishes asymptotically at the true parameter, its Hessian converges to a negative-definite matrix, and the empirical process determines the limiting distribution of the score.

Copulas provide a natural application because their defining objects are CDFs on a bounded domain. A copula separates dependence from marginal distributions; see \citet{nelsen2006}. Estimation based on an empirical copula is invariant under strictly increasing transformations of each margin. It is also important to distinguish an empirical distribution of observations from a copula with known margins from an empirical copula constructed using estimated margins. Their limiting processes differ, and the distinction affects the asymptotic variance of a fitted dependence parameter \citep{fermanian2004,tsukahara2005,segers2012}.

The closest related work is \citet{eguchikato2026}, who introduce minimum copula divergence estimators based on $\alpha$-, $\beta$- and $\gamma$-power losses and study their robustness. Our copula criterion belongs to the same broad class of distribution-function-based contrasts. In particular, its pointwise loss is a positive multiple of their $\beta$-power loss, apart from a term independent of the parameter. The main distinction is the integration measure. We integrate with respect to Lebesgue measure on the whole unit cube, whereas their criterion integrates with respect to the baseline copula measure and its empirical implementation is evaluated at the pseudo-observations. These choices generally yield different estimators and different covariance matrices. Our general cumulative formulation also applies beyond copulas. We therefore present the work as a study of a fixed-measure cumulative power criterion and its copula specialization, rather than as the introduction of an unrelated family of copula divergences.

Section~\ref{sec:criterion} defines the cumulative criterion and explains its domain of integration. Section~\ref{sec:asymptotic} states the consistency and asymptotic normality results. Section~\ref{sec:copulas} specializes them to copulas and gives a precise comparison with \citet{eguchikato2026}. Section~\ref{sec:simulations} describes the simulation design and the interpretation of its comparisons. Section~\ref{sec:conclusion} concludes. The proofs are collected in Appendix~\ref{sec:proofs}.

\section{The cumulative power criterion}\label{sec:criterion}

\subsection{From density power divergence to cumulative power divergence}

For two densities $g$ and $f_{\tht}$ on $\R^m$, the density power divergence of \citet{basu1998} is
\begin{equation}\label{c000}
d_a(g,f_{\tht})=\int_{\R^m}\left\{f_{\tht}^{1+a}
-\left(1+\frac1a\right)f_{\tht}^{a}g+\frac1a g^{1+a}\right\}\dd\boldsymbol{x},\qquad a>0.
\end{equation}
As $a\downarrow0$, it reduces to the Kullback--Leibler divergence, subject to the relevant integrability assumptions. At $a=1$ it is the squared $L_2$ discrepancy between the densities. Minimizing its empirical criterion at $a=0$ recovers maximum likelihood. That last correspondence concerns the density-based construction; the cumulative criterion below has a different limit and does not in general recover likelihood estimation.

Let $\boldsymbol{X}_1,\ldots,\boldsymbol{X}_n$ be independent and identically distributed random vectors in $\R^m$ with CDF $G$. Let $\{F_{\tht}:\tht\in\Theta\subseteq\R^p\}$ be a parametric family of CDFs. No density is assumed. Fix an integration domain $\D\subseteq\R^m$ and a finite, nonzero Borel measure $\mu$ on $\D$, independent of $\tht$ and of the observations. The main example is Lebesgue measure on a fixed bounded domain; for copulas, $\D=[0,1]^m$.
For $a>0$, define
\begin{equation}\label{c00}
\mathcal{C}d_{a,\mu}(G,F_{\tht})
=\int_{\D}\left\{F_{\tht}(\boldsymbol{x})^{1+a}
-\left(1+\frac1a\right)F_{\tht}(\boldsymbol{x})^aG(\boldsymbol{x})
+\frac1aG(\boldsymbol{x})^{1+a}\right\}\dd\mu(\boldsymbol{x}).
\end{equation}
We write $\mathcal{C}d_a$ when the measure is understood. This is the cumulative power construction of \citet{zografos2023}, with the integration measure stated explicitly. It is nonnegative because its pointwise integrand is $1/a$ times the Bregman remainder generated by the strictly convex function $T(t)=t^{1+a}$. It is zero precisely when $G=F_{\tht}$ $\mu$-almost everywhere. Throughout, identifiability is understood relative to this equality on $\D$. In particular, the chosen measure must distinguish the parameter values to be estimated. For Lebesgue measure on the unit cube, equality almost everywhere of two copulas implies equality everywhere by continuity.

The finite-measure convention ensures that the terms in~\eqref{c00} are finite. It is essential when the parameter-independent term is removed from the objective. On an unbounded domain equipped with Lebesgue measure, $\int F_{\tht}^{1+a}\dd\boldsymbol{x}$ is usually infinite, even if the combined divergence integrand is integrable. Thus such an unbounded-domain divergence cannot automatically be decomposed into the separate terms used below. A fixed integrable weight or a fixed bounded domain avoids subtracting infinite quantities.

At $a=1$, \eqref{c00} becomes
$
\mathcal{C}d_{1,\mu}(G,F_{\tht})=\int_{\D}\{G(\boldsymbol{x})-F_{\tht}(\boldsymbol{x})\}^2\dd\mu(\boldsymbol{x}).
$
The limiting pointwise loss as $a\downarrow0$ yields the cumulative Kullback--Leibler criterion
\begin{equation}\label{c01}
\mathcal{C}KL_{\mu}(G,F_{\tht})
=\int_{\D}\left\{G(\boldsymbol{x})\log\frac{G(\boldsymbol{x})}{F_{\tht}(\boldsymbol{x})}
+F_{\tht}(\boldsymbol{x})-G(\boldsymbol{x})\right\}\dd\mu(\boldsymbol{x}),
\end{equation}
when the limiting loss is integrable and the limit may be passed through the integral. We use the conventions $0\log(0/t)=0$ and an infinite loss if $G>0$ but $F_{\tht}=0$. Cumulative Kullback--Leibler estimation in the univariate setting is studied by \citet{mehraliasadi2021}. Our asymptotic results concern fixed $a>0$; the case $a=0$ is included as a numerical benchmark in Section~\ref{sec:simulations}.

\subsection{Empirical criterion and estimating equations}

Define the multivariate empirical CDF by
\begin{equation}\label{c1}
G_n(\boldsymbol{x})=\frac1n\sum_{i=1}^{n}\ind(\boldsymbol{X}_i\leq\boldsymbol{x}),
\end{equation}
where the inequality is coordinatewise. The minimum cumulative power divergence estimator (MCumPDE) is a minimizer of $\mathcal{C}d_{a,\mu}(G_n,F_{\tht})$. Since the term involving $G_n^{1+a}$ is independent of $\tht$, it can equivalently be defined by
\begin{equation}\label{c2}
\whth_{c,n}^{a}\in\underset{\tht\in\Theta}{\arg\min}\,V_n^a(\tht),
\end{equation}
where
\begin{equation}\label{c3}
V_n^a(\tht)=\int_{\D}\left\{F_{\tht}(\boldsymbol{x})^{1+a}
-\left(1+\frac1a\right)F_{\tht}(\boldsymbol{x})^aG_n(\boldsymbol{x})\right\}\dd\mu(\boldsymbol{x}).
\end{equation}
For the Taylor argument it is convenient to maximize
\begin{equation}\label{c3-1}
H_n^a(\tht)=-V_n^a(\tht),
\end{equation}
that is,
\begin{equation}\label{c3-2}
H_n^a(\tht)=\int_{\D}\left\{\left(1+\frac1a\right)F_{\tht}(\boldsymbol{x})^aG_n(\boldsymbol{x})
-F_{\tht}(\boldsymbol{x})^{1+a}\right\}\dd\mu(\boldsymbol{x}).
\end{equation}
Whenever differentiation under the integral is justified, an interior solution satisfies
\begin{equation}\label{c5}
\nabla_{\tht}H_n^a(\tht)
=(a+1)\int_{\D}w_a(\boldsymbol{x};\tht)\{G_n(\boldsymbol{x})-F_{\tht}(\boldsymbol{x})\}\dd\mu(\boldsymbol{x})=\boldsymbol{0}_p,
\end{equation}
with
\begin{equation}\label{eq:weight}
w_a(\boldsymbol{x};\tht)=F_{\tht}(\boldsymbol{x})^{a-1}\nabla_{\tht}F_{\tht}(\boldsymbol{x})
=F_{\tht}(\boldsymbol{x})^a\nabla_{\tht}\log F_{\tht}(\boldsymbol{x}).
\end{equation}
The first expression can be used wherever the logarithmic representation is inconvenient; zero-valued boundary points are handled by their limiting integrand when it exists. The parameter $a$ changes the weights of discrepancies at different cumulative probability levels. Its practical effect depends jointly on the model and on $\mu$.

\section{Consistency and asymptotic normality}\label{sec:asymptotic}

Suppose that $G=F_{\thzero}$ for an interior point $\thzero\in\Theta$, and fix $a>0$. We use the following local regularity conditions. They make explicit the domination needed in the Taylor arguments of the original cumulative criterion.

\begin{description}
\item[\textbf{(C1)}] The integration domain and measure are fixed and the family is identifiable relative to $\mu$. The true parameter $\thzero$ is an interior point of $\Theta$.
\item[\textbf{(C2)}] There is an open neighborhood $\Theta^*$ of $\thzero$ such that $F_{\tht}(\boldsymbol{x})>0$ for $\mu$-almost every $\boldsymbol{x}$ and all $\tht\in\Theta^*$, and $F_{\tht}(\boldsymbol{x})$ has continuous parameter derivatives through order three there.
\item[\textbf{(C3)}] The first and second parameter derivatives of
\[
h_a(\boldsymbol{x},z;\tht)=\left(1+\frac1a\right)F_{\tht}(\boldsymbol{x})^a z-F_{\tht}(\boldsymbol{x})^{1+a}
\]
are bounded by $\mu$-integrable envelopes, uniformly over $z\in[0,1]$ and $\tht$ in a closed neighborhood of $\thzero$. Differentiation under the integral is valid through order three, and $\int_{\D}\|w_a(\boldsymbol{x};\thzero)\|\dd\mu(\boldsymbol{x})<\infty$.
\item[\textbf{(C4)}] The finite $p\times p$ matrix $I_a(\thzero)$ with elements
\begin{equation}\label{c4}
I_{jk}(\tht)=\int_{\D}F_{\tht}(\boldsymbol{x})^{1+a}
\frac{\partial\log F_{\tht}(\boldsymbol{x})}{\partial\theta_j}
\frac{\partial\log F_{\tht}(\boldsymbol{x})}{\partial\theta_k}\dd\mu(\boldsymbol{x}),\quad j,k=1,\ldots,p,
\end{equation}
is positive definite. The subscript $a$ records its dependence on the fixed tuning parameter; it is not the likelihood Fisher information.
\item[\textbf{(C5)}] In a fixed closed neighborhood of $\thzero$, there are finite constants $M_{jk\ell}$, independent of $n$, such that
$\sup_n\sup_{\tht}\left|\frac{\partial^3H_n^a(\tht)}{\partial\theta_j\partial\theta_k\partial\theta_\ell}\right|\leq M_{jk\ell}$, $j,k,\ell=1,\ldots,p$.
%
A sufficient condition is a uniformly integrable envelope for the corresponding third derivatives of $h_a(\boldsymbol{x},z;\tht)$ over $z\in[0,1]$.
\end{description}

The domination assumptions are imposed on the weighted derivatives that occur in the criterion. Bounded derivatives of $F_{\tht}$ alone would not ensure these conditions when $a<1$ and $F_{\tht}$ approaches zero.

\begin{theorem}[Consistency of local solutions]\label{thm:consistency}
Under \textbf{(C1)--(C5)}, the estimating equations~\eqref{c5} admit a sequence of local maximizing solutions $\whth_{c,n}^a$, with probability tending to one, such that
\[
\whth_{c,n}^a\overset{P}{\longrightarrow}\thzero.
\]
Equivalently, for every $\zeta>0$, $P(\|\whth_{c,n}^a-\thzero\|<\zeta)\longrightarrow1$.
\end{theorem}

The proof in Appendix~\ref{proof:consistency} follows the local Taylor expansion approach of \citet{lehmanncasella1998}, \citet{basu2011} and \citet{castilla2022}. It establishes the existence of a consistent local solution, rather than consistency of every root of the estimating equations. To identify the global minimizer in~\eqref{c2} with a consistent solution, one can additionally restrict optimization to a compact parameter set containing $\thzero$ in its interior and assume pointwise continuity of $F_{\tht}$ on that set. With the finite measure and identifiability above, uniform convergence of $G_n$, continuity of the criterion and the usual argmin consistency theorem then give consistency of its global minimizers \citep[Theorem~5.7]{vandervaart1998}. The local normality result below applies to any consistent interior maximizing solution.
Write $F_0=F_{\thzero}$ and define
\begin{equation}\label{eq:Kcdf}
K_0(\boldsymbol{x},\boldsymbol{y})
=F_0(\boldsymbol{x}\wedge\boldsymbol{y})-F_0(\boldsymbol{x})F_0(\boldsymbol{y}),
\end{equation}
where $\boldsymbol{x}\wedge\boldsymbol{y}$ is the coordinatewise minimum. Set
\begin{equation}\label{eq:Bcdf}
B_a(\thzero)=\int_{\D}\int_{\D}K_0(\boldsymbol{x},\boldsymbol{y})
w_a(\boldsymbol{x};\thzero)w_a(\boldsymbol{y};\thzero)^\top \dd\mu(\boldsymbol{x})\dd\mu(\boldsymbol{y}).
\end{equation}
The integrability in \textbf{(C3)} makes this covariance finite, since $K_0$ is bounded. The product is the outer product of two column vectors.

\begin{theorem}[Asymptotic normality]\label{thm:normality}
Under \textbf{(C1)--(C5)}, every consistent sequence of interior solutions of~\eqref{c5} satisfies
\begin{equation}\label{eq:normal}
\sqrt{n}(\whth_{c,n}^a-\thzero)\overset{\mathcal{L}}{\longrightarrow}
N_p(\boldsymbol{0}_p,\Sigma_a),\qquad
\Sigma_a=I_a(\thzero)^{-1}B_a(\thzero)I_a(\thzero)^{-1}.
\end{equation}
\end{theorem}

The proof is in Appendix~\ref{proof:normality}. It uses the weak convergence of the multivariate empirical process indexed by lower orthants \citep{vandervaartwellner1996}. Both the Hessian and the score have a factor $a+1$, which cancels in~\eqref{eq:normal}. The covariance depends on $a$ through $I_a$ and $B_a$.

\section{Copula estimation and related divergences}\label{sec:copulas}

\subsection{The fixed Lebesgue-measure criterion}

Let $\boldsymbol{X}_i=(X_{i1},\ldots,X_{im})^\top$, $i = 1, \ldots, n$ have continuous marginal CDFs $F_1,\ldots,F_m$ and copula $C_0=C_{\thzero}$, so that
\[
F_0(\boldsymbol{x})=C_0\{F_1(x_1),\ldots,F_m(x_m)\}.
\]
If the margins are known, the transformed observations $\boldsymbol{U}_i=\{F_1(X_{i1}),\ldots,F_m(X_{im})\}^\top $ are independent observations from $C_0$. Their empirical CDF is
$
C_n^*(\boldsymbol{u})=\frac1n\sum_{i=1}^n\ind(\boldsymbol{U}_i\leq\boldsymbol{u}).
$
When the margins are unknown, use the pseudo-observations $\widehat U_{ij}=R_{ij}/(n+1)$, where $R_{ij}$ is the rank of $X_{ij}$ in the $j$th sample, and define
\begin{equation}\label{eq:rankcopula}
\widehat C_n(\boldsymbol{u})=\frac1n\sum_{i=1}^n\ind(\widehat{\boldsymbol{U}}_i\leq\boldsymbol{u}).
\end{equation}
We assume that ties have probability zero. The rank scaling keeps observations inside the open unit cube; the familiar asymptotically equivalent conventions give the same limiting process.
With $\D=[0,1]^m$ and $\mu$ equal to Lebesgue measure, define
\begin{equation}\label{r1}
\mathcal{C}d_a(C,C_{\tht})=
\int_{[0,1]^m}\left\{C_{\tht}(\boldsymbol{u})^{1+a}
-\left(1+\frac1a\right)C_{\tht}(\boldsymbol{u})^aC(\boldsymbol{u})
+\frac1aC(\boldsymbol{u})^{1+a}\right\}\dd\boldsymbol{u}.
\end{equation}
This cumulative divergence for copulas is discussed by \citet{zografos2024}. Its rank-based minimum copula power divergence estimator (MCopPDE) is
\begin{equation}\label{r2}
\whth_{cop,n}^a\in\underset{\tht\in\Theta}{\arg\min}
\int_{[0,1]^m}\left\{C_{\tht}(\boldsymbol{u})^{1+a}
-\left(1+\frac1a\right)C_{\tht}(\boldsymbol{u})^a\widehat C_n(\boldsymbol{u})\right\}\dd\boldsymbol{u}.
\end{equation}
Replacing $\widehat C_n$ by $C_n^*$ defines the corresponding known-margin estimator. Equation~\eqref{r2} is a minimization problem; equivalently, the negative of its displayed objective is maximized. At $a=1$ it minimizes the integrated squared difference between the empirical and model copulas with respect to Lebesgue measure.

The copula criterion is specified directly on the uniform scale. In general it is not obtained by simply changing variables in a Lebesgue-measure cumulative criterion on the original data scale: such a change of variables would introduce marginal weights. This distinction is part of the choice of integration measure.

\subsection{Known margins and estimated margins}

Put
\begin{align}\label{eq:Icop}
I_{a,cop}(\thzero)&=\int_{[0,1]^m}C_0(\boldsymbol{u})^{1+a}
s_0(\boldsymbol{u})s_0(\boldsymbol{u})^\top \dd\boldsymbol{u},\\
s_0(\boldsymbol{u})&=\nabla_{\tht}\log C_{\tht}(\boldsymbol{u})\big|_{\tht=\thzero},\qquad
w_{a,cop}(\boldsymbol{u})=C_0(\boldsymbol{u})^a s_0(\boldsymbol{u}).\notag
\end{align}
The derivatives in~\eqref{eq:Icop} are with respect to the model parameters, not to the coordinates of $\boldsymbol{u}$.

For known margins, the limiting process is the $C_0$-Brownian bridge $\mathbb{B}_{C_0}$, with covariance
\begin{equation}\label{eq:Kstar}
K^*(\boldsymbol{u},\boldsymbol{v})
=C_0(\boldsymbol{u}\wedge\boldsymbol{v})-C_0(\boldsymbol{u})C_0(\boldsymbol{v}).
\end{equation}
For ranks, write $\boldsymbol{e}_j(t)=(1,\ldots,1,t,1,\ldots,1)^\top $, with $t$ in position $j$. Under the empirical copula regularity conditions of \citet{segers2012},
\begin{equation}\label{eq:copprocess}
\sqrt n(\widehat C_n-C_0)\Longrightarrow\mathbb{C}_{C_0},\qquad
\mathbb{C}_{C_0}(\boldsymbol{u})=\mathbb{B}_{C_0}(\boldsymbol{u})
-\sum_{j=1}^m\dot C_{0,j}(\boldsymbol{u})\mathbb{B}_{C_0}\{\boldsymbol{e}_j(u_j)\}.
\end{equation}
Here $\dot C_{0,j}$ is a derivative with respect to the $j$th copula coordinate. Specifically, it is sufficient that these partial derivatives exist and are continuous on $\{\boldsymbol{u}:0<u_j<1\}$. The boundary values may be defined by the usual limiting conventions. Let
\begin{equation}\label{eq:Krank}
K^{rank}(\boldsymbol{u},\boldsymbol{v})
=\operatorname{Cov}\{\mathbb{C}_{C_0}(\boldsymbol{u}),\mathbb{C}_{C_0}(\boldsymbol{v})\}.
\end{equation}
Thus $K^{rank}$ is obtained by expanding~\eqref{eq:copprocess} using~\eqref{eq:Kstar}; it includes the marginal correction terms.

\begin{proposition}[Copula specialization]\label{prop:copula}
Assume the copula versions of \textbf{(C1)--(C5)} on the unit cube. The known-margin criterion admits a consistent sequence of local solutions, and its consistent interior solutions satisfy
\[
\sqrt n(\whth_{cop,n}^{a,*}-\thzero)\overset{\mathcal{L}}{\longrightarrow}
N_p(\boldsymbol{0}_p,I_{a,cop}^{-1}B_{a,cop}^*I_{a,cop}^{-1}),
\]
where, with the matrices evaluated at $\thzero$,
\begin{equation}\label{eq:Bcop}
B_{a,cop}^*=
\int_{[0,1]^m}\int_{[0,1]^m}K^*(\boldsymbol{u},\boldsymbol{v})
w_{a,cop}(\boldsymbol{u})w_{a,cop}(\boldsymbol{v})^\top \dd\boldsymbol{u}\dd\boldsymbol{v}.
\end{equation}
If, in addition, the empirical copula conditions in~\eqref{eq:copprocess} hold, the same conclusions apply to the rank-based criterion, with $B_{a,cop}^*$ replaced by $B_{a,cop}^{rank}$, defined by replacing $K^*$ with $K^{rank}$ in~\eqref{eq:Bcop}.
\end{proposition}

This proposition uses the same Taylor expansions as Theorems~\ref{thm:consistency} and~\ref{thm:normality}. For known margins it is their direct specialization to $F_{\tht}=C_{\tht}$. For estimated margins, uniform consistency of $\widehat C_n$ gives the identical first-derivative and Hessian limits; its weak limit~\eqref{eq:copprocess} replaces the empirical-distribution limit in the last step of the normality argument. The integrability of $w_{a,cop}$ makes this integration a continuous linear map. No additional expansion is needed. These covariance formulas apply under correct specification; the contaminated and misspecified models in the simulations are assessed through finite-sample errors rather than through these normal limits.

\subsection{Relationship with Eguchi and Kato}

For $a>0$ define the pointwise $\beta$-power loss
$
\ell_a^{\beta}(c_0,c_1)=\frac{c_1^{1+a}}{1+a}-\frac{c_0c_1^a}{a}.
$
The loss in~\eqref{r1} satisfies the exact identity
\begin{equation}\label{eq:betaidentity}
c_1^{1+a}-\left(1+\frac1a\right)c_0c_1^a+\frac1a c_0^{1+a}
=(1+a)\ell_a^{\beta}(c_0,c_1)+\frac1a c_0^{1+a}.
\end{equation}
Consequently, replacing the name of the loss or multiplying it by a positive constant does not produce a different estimator if the same integration measure is used.

The $\beta$-copula divergence of \citet{eguchikato2026} is
\begin{equation}\label{eq:eguchibeta}
D_a^{EK}(C,C_{\tht})=
\int_{[0,1]^m}\left\{\frac{C_{\tht}^{1+a}}{1+a}
-\frac{CC_{\tht}^{a}}{a}+\frac{C^{1+a}}{a(1+a)}\right\}\dd C(\boldsymbol{u}).
\end{equation}
Their empirical criterion evaluates the loss at the pseudo-observations:
\begin{equation}\label{eq:EKempirical}
\frac1n\sum_{i=1}^n\ell_a^{\beta}\{\widehat C_n(\widehat{\boldsymbol{U}}_i),C_{\tht}(\widehat{\boldsymbol{U}}_i)\}.
\end{equation}
Our criterion instead integrates the same pointwise loss over uniformly weighted locations in the cube. Table~\ref{tab:comparison} summarizes the distinction.

\begin{table}[tbp]\centering\small
\caption{The role of the integration measure for the same power loss.}\label{tab:comparison}
\begin{tabular}{@{}p{0.23\linewidth}p{0.35\linewidth}p{0.35\linewidth}@{}}\toprule
 & Present copula criterion & Eguchi--Kato $\beta$ criterion\\\midrule
Population measure & Lebesgue measure $\dd\boldsymbol{u}$ & Baseline copula measure $\dd C(\boldsymbol{u})$\\
Empirical evaluation & Quadrature over the unit cube & Average over pseudo-observations\\
Pointwise loss & $(1+a)\ell_a^{\beta}$, up to a parameter-independent term & $\ell_a^{\beta}$, up to a parameter-independent term\\
Correct-model Hessian up to positive scaling& $\int C_0^{1+a}s_0s_0^\top \dd\boldsymbol{u}$ & $\int C_0^{1+a}s_0s_0^\top \dd C_0(\boldsymbol{u})$\\
General setting & CDFs and copulas with a fixed finite measure & Copula divergence estimation\\\bottomrule
\end{tabular}
\end{table}

The difference is in spatial weighting. Lebesgue integration assesses discrepancies throughout the cube; integration against $C$ emphasizes locations according to the copula distribution. This difference remains when densities exist and is particularly transparent at $a=1$, where the two procedures use different integrated squared-error criteria. Under misspecification, they also minimize different population objectives and can have different pseudo-true parameters. Neither criterion is universally preferable without specifying the statistical objective and the type of misspecification.

The two approaches share their use of copula functions and empirical copulas. They also share invariance to strictly increasing marginal transformations when constructed from ranks. Therefore that invariance, density-free evaluation and minimum copula divergence estimation are not claimed here as new properties relative to \citet{eguchikato2026}. Their power-boundedness and robustness analysis provides useful complementary results. We concentrate on the cumulative fixed-measure criterion, its asymptotic covariance and the numerical consequences of the weighting choice; we do not assert a general robustness dominance over their estimators.

\section{Simulation study}\label{sec:simulations}

\subsection{Design and target of evaluation}

We consider Gumbel, Clayton and Frank copulas with Kendall's tau
$\tau_0\in\{0.1,0.3,0.5,0.9\}$, sample sizes $n\in\{100,1000\}$ and
contamination fractions $\varepsilon\in\{0,0.1,0.5\}$. Each generating
family is fitted by all three families, giving 216 design cells with
200 replications each. The parameters are
$\theta_G(\tau)=1/(1-\tau)$ and $\theta_C(\tau)=2\tau/(1-\tau)$;
the Frank parameter is obtained by numerical inversion. These families
differ in tail dependence despite having equal Kendall tau
\citep{nelsen2006}.

Each sample pools $n_1=\operatorname{round}(n\varepsilon)$ independent
uniform pairs and $n-n_1$ draws from $C_{\theta_0}$, corresponding to
\begin{equation}\label{eq:mixture}
C_\varepsilon=(1-\varepsilon)C_{\theta_0}+\varepsilon\Pi,
\qquad \Pi(u,v)=uv.
\end{equation}
Within each replicate at fixed $n$, $\tau_0$ and generating family,
clean and independent draws are reused across contamination fractions;
all methods and fitted families receive the same pooled sample.

The evaluation target is the uncontaminated dependence:
\begin{equation}\label{eq:abserror}
E=|\widehat\tau-\tau_0|.
\end{equation}
This differs from estimating the mixture's Kendall tau, which is
generally not $(1-\varepsilon)\tau_0$. Under family misspecification,
the divergence minimizer may also have Kendall tau different from
$\tau_0$. Persistent errors in these settings therefore do not
contradict consistency under correct specification.

\subsection{Estimators and computation}

We evaluate MCopPDE at $a\in\{0.1,0.5,1\}$ and include its cumulative
Kullback--Leibler limit $a=0$ as a numerical benchmark. Fits use the
unsmoothed empirical copula based on marginal ranks
$\widehat U_{ij}=R_{ij}/(n+1)$. The criterion approximates Lebesgue
integration by summing losses over a uniform $1001\times1001$ grid
on $[0,1]^2$, with spacing $h=0.001$.

The comparison estimators of \citet{eguchikato2026}, denoted EK--$\alpha$,
EK--$\beta$ and EK--$\gamma$, evaluate their losses at the sample
pseudo-observations and integrate with respect to the empirical copula.
We use $q\in\{0.1,0.5,1\}$ for EK--$\beta$ and EK--$\gamma$, and
$q\in\{0.1,0.5\}$ for EK--$\alpha$. At $q=a$, the pointwise $\beta$
and cumulative power losses are equivalent up to scaling and
parameter-independent terms; their integration measures differ here.
All divergence criteria are minimized over
$\tau\in[10^{-5},1-10^{-5}]$, transforming trial values to family
parameters. A full-interval search is supplemented by a nine-point
scan and searches in the two best neighboring brackets, retaining
the smallest evaluated criterion.

We also record oracle maximum likelihood fits to the generated
uniforms, maximum pseudo-likelihood fits to the ranks, and untruncated
sample Kendall tau. Likelihood fits use the package's parameter
domains, whereas divergence fits allow only positive dependence;
this distinction can affect weak-dependence comparisons.

\subsection{Performance summaries and numerical diagnostics}

We report MAE, the average of~\eqref{eq:abserror}, and signed bias,
the average of $\widehat\tau-\tau_0$. For $B$ available replications,
the MCSE of MAE is $s_E/\sqrt B$, where $s_E$ is the sample standard
deviation of absolute errors. Method differences use paired MCSEs.

All divergence fits returned finite estimates without recorded errors.
Two oracle likelihood fits and one pseudo-likelihood fit failed,
leaving 199 values in each affected cell. Likelihood convergence
diagnostics were not saved, and many fitted-Clayton likelihood
estimates numerically coincide with sample Kendall tau. Comparisons
therefore focus on divergence estimators and sample Kendall tau.

\begin{figure}[!htbp]\centering
\includegraphics[width=\linewidth]{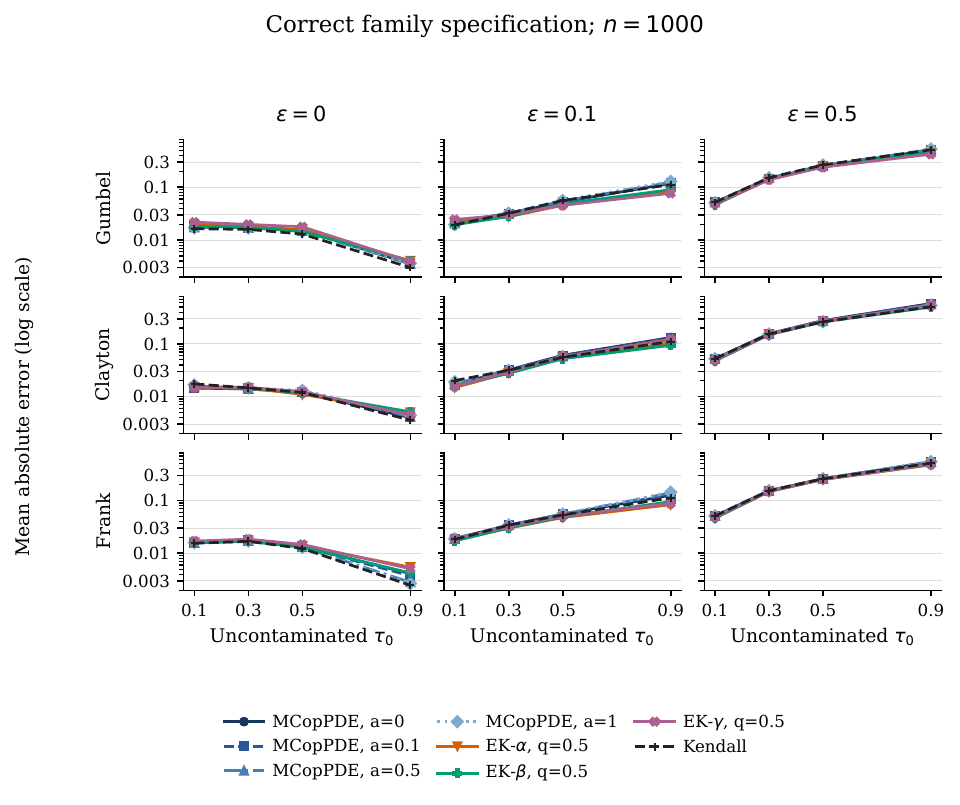}
\caption{Correctly specified fits at $n=1000$. Curves show mean absolute error relative to the uncontaminated $\tau_0$, based on 200 replications per point. Rows identify the generating and fitted family; columns identify contamination. The vertical scale is logarithmic and common to all panels. All four MCopPDE powers are shown, together with the three EK estimators at the common intermediate exponent $q=0.5$ and sample Kendall tau. Results for all EK powers and $n=100$ are reported in Table~\ref{tab:simulation-overview} and the supplement.}\label{fig:sim-correct}
\end{figure}

Table~\ref{tab:simulation-overview} covers all evaluated powers and
both sample sizes, equally averaging the 12 correctly specified cells
at each $(n,\varepsilon)$. For independent cell means, its MCSE is
$\{\sum_{j=1}^{12}\operatorname{MCSE}_j^2\}^{1/2}/12$.
Figure~\ref{fig:sim-correct} gives family-specific results under correct
specification. Supplementary Tables~S1
and~S2 compare MCopPDE and EK--$\beta$
at $a=q=0.5$ for all family combinations, dependence levels and
contamination fractions at $n=1000$; Supplementary Tables~S3--S4
give the corresponding $n=100$ results. Supplementary Tables~S5--S22
report MAEs, signed biases and MCSEs for all 216 cells and evaluated
divergence powers, together with sample Kendall tau.
\vspace*{-5mm}
\begin{table}[!htbp]
\centering\small
\caption{Mean absolute errors under correct family specification, equally averaged over three families and four values of $\tau_0$. Each entry averages 12 cell means, each based on 200 replications. Parentheses give the Monte Carlo standard error of this average in units of $10^{-4}$, rounded to the nearest integer.}
\label{tab:simulation-overview}
\setlength{\tabcolsep}{4pt}
\begin{tabular}{@{}lc rr rr rr@{}}\toprule
& & \multicolumn{2}{c}{$\varepsilon=0$} & \multicolumn{2}{c}{$\varepsilon=0.1$} & \multicolumn{2}{c}{$\varepsilon=0.5$}\\
\cmidrule(lr){3-4}\cmidrule(lr){5-6}\cmidrule(lr){7-8}
Estimator & Power & $n=100$ & $n=1000$ & $n=100$ & $n=1000$ & $n=100$ & $n=1000$\\\midrule
MCopPDE & 0 & 0.0437 (7) & 0.0126 (2) & 0.0866 (9) & 0.0584 (3) & 0.2642 (12) & 0.2522 (4)\\
MCopPDE & 0.1 & 0.0429 (7) & 0.0125 (2) & 0.0856 (9) & 0.0585 (3) & 0.2630 (12) & 0.2523 (4)\\
MCopPDE & 0.5 & 0.0411 (7) & 0.0123 (2) & 0.0824 (9) & 0.0588 (3) & 0.2589 (12) & 0.2520 (4)\\
MCopPDE & 1 & 0.0406 (7) & 0.0125 (2) & 0.0800 (9) & 0.0587 (3) & 0.2544 (12) & 0.2509 (4)\\
\addlinespace[2pt]
EK--$\alpha$ & 0.1 & 0.0581 (8) & 0.0131 (2) & 0.0556 (8) & 0.0463 (3) & 0.2025 (12) & 0.2308 (4)\\
EK--$\alpha$ & 0.5 & 0.0567 (8) & 0.0130 (2) & 0.0556 (8) & 0.0465 (3) & 0.2054 (12) & 0.2318 (4)\\
\addlinespace[2pt]
EK--$\beta$ & 0.1 & 0.0571 (8) & 0.0129 (2) & 0.0552 (8) & 0.0465 (3) & 0.2037 (12) & 0.2315 (4)\\
EK--$\beta$ & 0.5 & 0.0551 (8) & 0.0126 (2) & 0.0549 (8) & 0.0473 (3) & 0.2075 (11) & 0.2332 (4)\\
EK--$\beta$ & 1 & 0.0566 (8) & 0.0128 (2) & 0.0552 (8) & 0.0475 (3) & 0.2072 (11) & 0.2331 (4)\\
\addlinespace[2pt]
EK--$\gamma$ & 0.1 & 0.0545 (9) & 0.0147 (3) & 0.0656 (9) & 0.0473 (4) & 0.2091 (14) & 0.2279 (5)\\
EK--$\gamma$ & 0.5 & 0.0481 (8) & 0.0138 (2) & 0.0655 (9) & 0.0494 (3) & 0.2203 (13) & 0.2332 (5)\\
EK--$\gamma$ & 1 & 0.0455 (8) & 0.0134 (2) & 0.0677 (9) & 0.0510 (3) & 0.2273 (12) & 0.2360 (4)\\
\addlinespace[2pt]
Kendall & -- & 0.0395 (7) & 0.0119 (2) & 0.0732 (9) & 0.0547 (3) & 0.2476 (12) & 0.2443 (4)\\
\bottomrule\end{tabular}
\end{table}
\vspace*{-10mm}

\subsection{Results}

Under correct specification without contamination, the cumulative estimators are competitive with the rank-based alternatives. For example, the average MAE of MCopPDE with $a=0.5$ decreases from 0.0411 at $n=100$ to 0.0123 at $n=1000$; the corresponding values are 0.0551 and 0.0126 for EK--$\beta$ with $q=0.5$, and 0.0395 and 0.0119 for Kendall tau. Thus the advantage over this EK--$\beta$ fit at $n=100$ becomes small at $n=1000$, while Kendall tau has slightly smaller errors in the averaged clean-data comparison. These averages do not imply a uniform ordering. At $n=100$ and $\tau_0=0.9$, correctly specified Frank fits have MAE 0.0326 at $a=0$, compared with 0.0119 at $a=1$ and 0.0086 for Kendall tau, illustrating a finite-sample cost of the cumulative limit at strong dependence (Supplementary Table~S13).

Independence contamination substantially increases errors relative to the uncontaminated target. At $n=1000$, the average MAE of MCopPDE with $a=0.5$ is 0.0588 for $\varepsilon=0.1$ and 0.2520 for $\varepsilon=0.5$, compared with 0.0473 and 0.2332 for the matched EK--$\beta$ estimator. The latter has smaller observed MAE in all 24 correctly specified contaminated cells at this sample size. In particular, at $\tau_0=0.9$ and $\varepsilon=0.1$, MCopPDE with $a=0.5$ has MAEs between 0.1191 and 0.1396 across the three families. The smaller errors on clean data therefore do not translate into superior recovery under this contamination mechanism. At $\varepsilon=0.5$, all displayed methods have substantial recovery errors. Increasing the sample size does not remove the discrepancy between the clean-component target and the contaminated population; some contaminated-cell errors even increase as finite-sample biases change.

The effect of the power parameter depends on the family and scenario. For correctly specified Clayton data with $n=1000$, $\tau_0=0.9$ and $\varepsilon=0.1$, increasing $a$ from 0 to 1 reduces MAE from 0.1319 to 0.1087. In the corresponding Frank case it increases MAE from 0.1265 to 0.1439. Hence larger powers do not provide a general improvement in recovery, and these results do not identify one uniformly preferable value of $a$ (Supplementary Tables~S18 and~S22).

\begin{figure}[!htbp]\centering
\includegraphics[width=\linewidth]{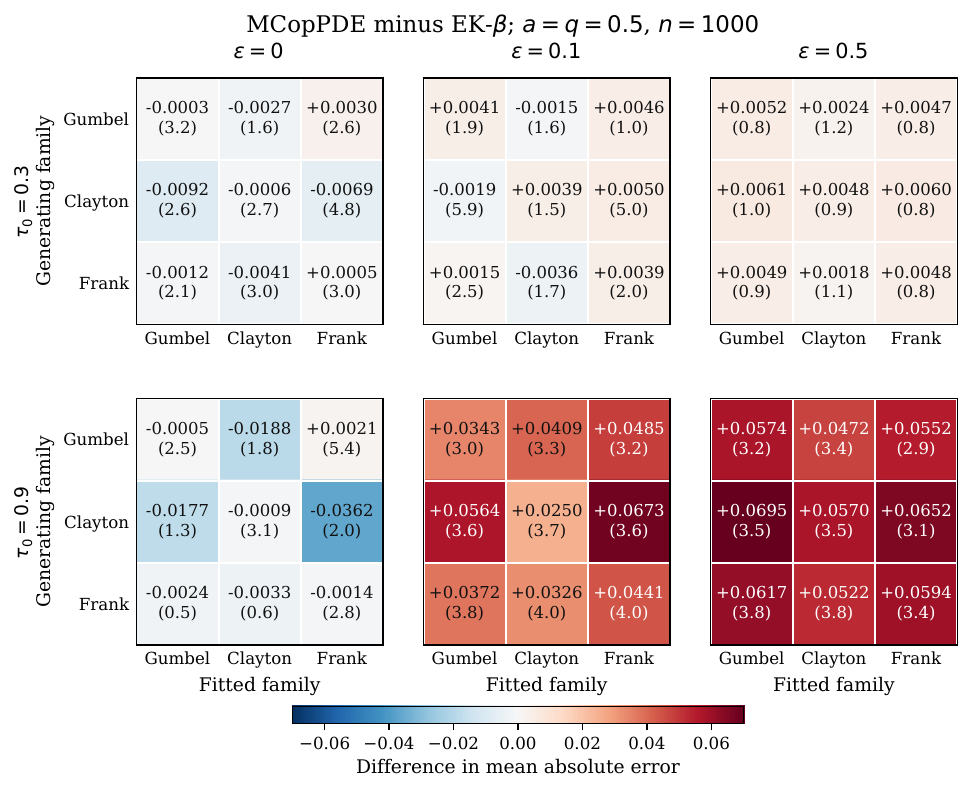}
\caption{Differences in mean absolute error, MCopPDE minus EK--$\beta$, for $a=q=0.5$ and $n=1000$. Each matrix has generating families in rows and fitted families in columns; diagonal entries provide correctly specified references. Blue (negative) entries favor MCopPDE and red (positive) entries favor EK--$\beta$. Each cell gives the MAE difference, followed in parentheses by its paired MCSE in units of $10^{-4}$. All cells use 200 paired replications. Errors refer to the uncontaminated $\tau_0$.}\label{fig:sim-weights}
\end{figure}

Figure~\ref{fig:sim-weights} compares MCopPDE and EK--$\beta$ at the common intermediate exponent $a=q=0.5$, for moderate and strong dependence. It includes all nine family combinations and all contamination fractions. The differences change sign across scenarios. For instance, with $n=1000$ and $\tau_0=0.9$, fitting Frank to clean Clayton data gives MAE 0.0100 for MCopPDE and 0.0463 for EK--$\beta$; at $\varepsilon=0.1$ the errors become 0.1334 and 0.0660 (Supplementary Tables~S2 and S19). Thus the integration measure materially affects the fitted dependence, and a favorable ordering on clean misspecified data can reverse under contamination.

Occasional decreases in error after adding contamination require particular care under misspecification. For example, fitting Gumbel to Clayton data at $n=1000$ and $\tau_0=0.3$, MCopPDE with $a=0.5$ has MAE 0.0372 without contamination and 0.0160 at $\varepsilon=0.1$. The corresponding biases are $+0.0365$ and $-0.0005$, indicating cancellation between the misspecification bias and the shift caused by contamination (see Supplementary Tables~S1 and  S17). Such a result does not establish general resistance to contamination.

Overall, the simulations show the numerical consequences of changing the integration measure and the limitations of clean-component recovery under independence contamination. MCopPDE performs reasonably on clean, correctly specified samples, but the present design provides no evidence for a general robustness advantage over the Eguchi--Kato estimators. The complete scenario results are retained in the supplement, including unfavorable cases and all tested powers.

\section{Concluding remarks}\label{sec:conclusion}

Minimum cumulative power divergence estimation fits multivariate
distribution functions without density smoothing. For a fixed positive
power and a fixed finite integration measure, the estimating equations
admit consistent, asymptotically normal local solutions under the stated
conditions. For copulas under correct specification, known and estimated
margins yield the same limiting Hessian but different score covariances.

At the same positive exponent, the pointwise loss is equivalent to the
$\beta$ contrast of \citet{eguchikato2026}, up to scaling and
parameter-independent terms. The contribution lies in the fixed-measure
formulation and its asymptotic analysis; our copula implementation uses
Lebesgue integration rather than empirical copula weighting.

The simulations show competitive performance on clean, correctly specified
data, but provide no evidence of a general robustness advantage under
independence contamination. No evaluated power is uniformly preferable.
Errors target the uncontaminated Kendall tau, so persistent errors under
contamination or misspecification do not contradict consistency under
correct specification. Improvements after contamination may reflect bias
cancellation, and accurate Kendall tau estimates need not imply accurate
recovery of the full copula.

Further work includes data-driven power selection, alternative fixed
weights, higher-dimensional integration, and evaluation beyond Kendall's tau.

\section*{Statements and Declarations}

\textbf{Funding:} This work began during the second author's December 2024 visit to
Technische Universit\"at Dresden, supported through a Greek--German
exchange programme by the Greek Ministry of Education, Religious
Affairs and Sports and the German Academic Exchange Service (DAAD).\\
\textbf{Competing interests:} The authors declare no competing interests.\\
\textbf{Data availability:} This study uses simulated data only. 
\begin{appendices}

\section{Proofs of the asymptotic results}\label{sec:proofs}
We give the full Taylor expansions, intermediate bounds and empirical-process calculations for the two asymptotic results. Throughout, integration is with respect to the fixed finite measure $\mu$ on $\D$, and the true distribution is $F_{\thzero}$.
\subsection{Proof of Theorem~\ref{thm:consistency}}\label{proof:consistency}

To start the proof, following the exposition in \citet[p.~304]{basu2011} and \citet[p.~4604]{castilla2022}, we will study the behavior of the
objective function $H_{n}^{a}(\boldsymbol{\theta })$, defined by (\ref{c3-2}), namely, 
\begin{equation}
H_{n}^{a}(\boldsymbol{\theta })=\dint\limits_{\D}\left\{ \left( 1+\frac{1}{a}\right) F_{\boldsymbol{\theta }}(
\boldsymbol{x})^{a}G_{n}(\boldsymbol{x})-F_{\boldsymbol{\theta }}(
\boldsymbol{x})^{1+a}\right\} d\mu(\boldsymbol{x}),\text{ }a>0,  \label{c3-3}
\end{equation}
on a neighbourhood of the true parameter value $\boldsymbol{\theta }_{0}$.
Let $Q_{\zeta }$ denote the closed Euclidean ball with center at $\boldsymbol{\theta }_{0}$
and radius $\zeta $. The interest is focused in showing that for any
sufficiently small $\zeta $, $H_{n}^{a}(\boldsymbol{\theta })<H_{n}^{a}(
\boldsymbol{\theta }_{0})$ with probability tending to 1 at all points $
\boldsymbol{\theta }$ on the surface of $Q_{\zeta }$. Following \citet[p.~4604]{castilla2022}, and adapting the argument to the present notation, this inequality ensures that the objective function $H_n^a(\boldsymbol{\theta})$ has a maximum in the interior of $Q_\zeta$. Indeed, continuity gives a maximizer on the closed ball, and the strict boundary inequality excludes its surface, since $\boldsymbol{\theta}_0$ belongs to the ball. Such a maximizer is a local maximum of the criterion. Since $H_n^a(\boldsymbol{\theta})$ is differentiable, the system of equations~\eqref{c5} must be satisfied at an interior local maximum. Therefore, for every sufficiently small $\zeta>0$, the system has, with probability tending to one, a solution $\widehat{\boldsymbol{\theta}}_{c,n}^a=\widehat{\boldsymbol{\theta}}_{c,n}^a(\zeta)$ within $Q_\zeta$, verifying
$
\underset{n\rightarrow \infty }{\lim }P\left( \left\Vert \widehat{
\boldsymbol{\theta }}_{c,n}^{a}-\boldsymbol{\theta }_{0}\right\Vert <\zeta
\right) =1.
$
Consider a Taylor series expansion of $H_{n}^{a}(\boldsymbol{\theta })$
around the true value $\boldsymbol{\theta }_{0}=(\theta _1^0,\ldots,\theta_p^0)^\top ,$
\begin{align}
H_{n}^{a}(\boldsymbol{\theta })-H_{n}^{a}(\boldsymbol{\theta }_{0})
&=\sum\limits_{j=1}^{p}\left( \frac{\partial H_{n}^{a}(\boldsymbol{\theta })
}{\partial \theta _{j}}\right) _{\boldsymbol{\theta =\theta }_{0}}(\theta
_{j}-\theta _{j}^{0})  \notag \\
&\quad+\frac{1}{2}\sum\limits_{j=1}^{p}\sum\limits_{k=1}^{p}\left( \frac{
\partial ^{2}H_{n}^{a}(\boldsymbol{\theta })}{\partial \theta _{j}\partial
\theta _{k}}\right) _{\boldsymbol{\theta =\theta }_{0}}(\theta _{j}-\theta
_{j}^{0})(\theta _{k}-\theta _{k}^{0})  \label{c6} \\
&\quad+\frac{1}{6}\sum\limits_{j=1}^{p}\sum\limits_{k=1}^{p}\sum\limits_{\ell
=1}^{p}\left( \frac{\partial ^{3}H_{n}^{a}(\boldsymbol{\theta })}{\partial
\theta _{j}\partial \theta _{k}\partial \theta _{\ell }}\right) _{
\boldsymbol{\theta =\theta }^{\ast }}(\theta _{j}-\theta _{j}^{0})(\theta
_{k}-\theta _{k}^{0})(\theta _{\ell }-\theta _{\ell }^{0})  \notag \\
&=L_{1}+L_{2}+L_{3},  \notag
\end{align}
where $\boldsymbol{\theta}^{\ast}$ lies on the segment joining $\boldsymbol{\theta}_0$ to $\boldsymbol{\theta}$ and hence in the ball of radius $\zeta$. In the sequel,
we will manipulate separately each term of the right-hand side of (\ref{c6}).

We start with $L_{1}$. Based on (\ref{c3-3}) and taking into account
condition (\textbf{C3}), 
\begin{align}
\frac{\partial H_{n}^{a}(\boldsymbol{\theta })}{\partial \theta _{j}}
&=\dint\limits_{\D}\left\{ \left( 1+\frac{1}{a}\right) aF_{\boldsymbol{\theta }}(
\boldsymbol{x})^{a-1}G_{n}(\boldsymbol{x})\frac{\partial F_{\boldsymbol{
\theta }}(\boldsymbol{x})}{\partial \theta _{j}}-(1+a)F_{\boldsymbol{\theta }
}(\boldsymbol{x})^{a}\frac{\partial F_{\boldsymbol{\theta }}(\boldsymbol{x})
}{\partial \theta _{j}}\right\} d\mu(\boldsymbol{x}),  \notag \\
&=(1+a)\dint\limits_{\D}\left\{ F_{\boldsymbol{\theta }}(\boldsymbol{x})^{a-1}\frac{\partial F_{
\boldsymbol{\theta }}(\boldsymbol{x})}{\partial \theta _{j}}[G_{n}(
\boldsymbol{x})-F_{\boldsymbol{\theta }}(\boldsymbol{x})]\right\} d\mu(\boldsymbol{x}).  \label{c7}
\end{align}
\noindent 
%
By the multivariate Glivenko--Cantelli theorem,
$\|G_n-F_0\|_\infty\to0$ almost surely. Together with \textbf{(C3)},
this yields:
\begin{equation}
\left( \frac{\partial H_{n}^{a}(\boldsymbol{\theta })}{\partial \theta _{j}}
\right) _{\boldsymbol{\theta =\theta }_{0}}\underset{n\rightarrow \infty }{
\overset{P}{\rightarrow }}0.  \label{c8}
\end{equation}
Hence, for any given $\zeta>0$, it follows, as in \citet[p.~431]{lehmann1983} or \citet[p.~464]{lehmanncasella1998}, that the event
$
\max_{1\leq j\leq p}\left|\left(\frac{\partial H_n^a(\boldsymbol{\theta})}{\partial\theta_j}\right)_{\boldsymbol{\theta}=\boldsymbol{\theta}_0}\right|<\zeta^2
$
has probability tending to one as $n\to\infty$. There are finitely many coordinates, so their convergence in probability in~\eqref{c8} gives the simultaneous bound. In particular, the bound concerns the absolute values of the derivatives; their signs are not restricted. All the bounds for $L_1$ below hold on this event, uniformly over the surface of $Q_\zeta$.
In this context,
$|L_1|
=\left|\sum_{j=1}^p
\left(\frac{\partial H_n^a(\boldsymbol{\theta})}{\partial\theta_j}\right)_{\boldsymbol{\theta}=\boldsymbol{\theta}_0}
(\theta_j-\theta_j^0)\right|
\leq\sum_{j=1}^p
\left|\left(\frac{\partial H_n^a(\boldsymbol{\theta})}{\partial\theta_j}\right)_{\boldsymbol{\theta}=\boldsymbol{\theta}_0}\right|
|\theta_j-\theta_j^0|
\leq\zeta\sum_{j=1}^p
\left|\left(\frac{\partial H_n^a(\boldsymbol{\theta})}{\partial\theta_j}\right)_{\boldsymbol{\theta}=\boldsymbol{\theta}_0}\right|
$
and taking into account $\left|\left(\partial H_n^a(\boldsymbol{\theta})/\partial\theta_j\right)_{\boldsymbol{\theta}=\boldsymbol{\theta}_0}\right|<\zeta^2$, for every $j$, with probability tending to $1$,
\begin{equation}
\left\vert L_{1}\right\vert \leq \zeta \sum\limits_{j=1}^{p}\zeta
^{2}=p\zeta ^{3}.  \label{c8a}
\end{equation}
We now continue with $L_{2}$. Based on (\ref{c3-3}), (\ref{c6}), (\ref{c7}
), after simple algebraic manipulations,
\begin{align}
\frac{\partial ^{2}H_{n}^{a}(\boldsymbol{\theta })}{\partial \theta
_{j}\partial \theta _{k}} &=(a+1)\dint\limits_{\D}\frac{\partial }{\partial \theta _{k}}\left\{ F_{\boldsymbol{\theta }}(
\boldsymbol{x})^{a-1}\frac{\partial }{\partial \theta _{j}}F_{\boldsymbol{
\theta }}(\boldsymbol{x})\right\} [G_{n}(\boldsymbol{x})-F_{\boldsymbol{
\theta }}(\boldsymbol{x})]d\mu(\boldsymbol{x})  \notag \\
&\quad-(a+1)\dint\limits_{\D}F_{\boldsymbol{\theta }}(\boldsymbol{x})^{a+1}\frac{\partial }{\partial
\theta _{j}}\log F_{\boldsymbol{\theta }}(\boldsymbol{x})\frac{\partial }{
\partial \theta _{k}}\log F_{\boldsymbol{\theta }}(\boldsymbol{x})d\mu(\boldsymbol{x}).  \label{c9}
\end{align}
Given that $G_{n}(\boldsymbol{x})\rightarrow F_{\boldsymbol{\theta }_{0}}(
\boldsymbol{x})$, with probability one, together with the domination in \textbf{(C3)}, (\ref{c9}) above leads to 
\begin{equation}
\left( \frac{\partial ^{2}H_{n}^{a}(\boldsymbol{\theta })}{\partial \theta
_{j}\partial \theta _{k}}\right) _{\boldsymbol{\theta =\theta }
_{0}}\rightarrow -(a+1)I_{jk}(\boldsymbol{\theta }_{0}),  \label{c9a}
\end{equation}
where $I_{jk}(\boldsymbol{\theta })$ are the $(j,k)$-elements of the positive
definite matrix of condition (\textbf{C4}). Following the steps of the proof in \citet[p.~305]{basu2011} and taking into account (\ref{c6}),
\begin{align}
2L_{2} &=\dsum\limits_{j=1}^{p}\dsum\limits_{k=1}^{p}\left( \frac{\partial
^{2}H_{n}^{a}(\boldsymbol{\theta })}{\partial \theta _{j}\partial \theta _{k}
}\right) _{\boldsymbol{\theta }=\boldsymbol{\theta }_{0}}(\theta _{j}-\theta
_{j}^{0})(\theta _{k}-\theta _{k}^{0})  \notag \\
&=\dsum\limits_{j=1}^{p}\dsum\limits_{k=1}^{p}\left( -(a+1)I_{jk}(
\boldsymbol{\theta }_{0})(\theta _{j}-\theta _{j}^{0})(\theta _{k}-\theta
_{k}^{0})\right)   \notag \\
&\quad+\dsum\limits_{j=1}^{p}\dsum\limits_{k=1}^{p}\left[ \left( \frac{\partial
^{2}H_{n}^{a}(\boldsymbol{\theta })}{\partial \theta _{j}\partial \theta _{k}
}\right) _{\boldsymbol{\theta }=\boldsymbol{\theta }_{0}}-\left[
-(a+1)I_{jk}(\boldsymbol{\theta }_{0})\right] \right] (\theta _{j}-\theta
_{j}^{0})(\theta _{k}-\theta _{k}^{0}).  \label{c10}
\end{align}
As in \citet[pp.~431--432]{lehmann1983} or \citet[p.~464]{lehmanncasella1998},
for the second term of (\ref{c10}), it follows from an argument similar to
that of $L_{1}$ that its absolute value is less than $p^{2}\zeta ^{3}$, $
\zeta >0,$ with probability tending to $1$. Closely following the exposition
in \citet[p.~432]{lehmann1983}, the first term in (\ref{c10})
is a negative (nonrandom) quadratic form in the variables $(\theta
_{j}-\theta _{j}^{0})$. By an orthogonal transformation the quadratic form $
\sum_{j=1}^{p}\sum_{k=1}^{p}\left( -\frac{1}{2}(a+1)I_{jk}(
\boldsymbol{\theta }_{0})(\theta _{j}-\theta _{j}^{0})(\theta _{k}-\theta
_{k}^{0})\right) $ can be reduced to diagonal form $\sum_{i=1}^{p}
\lambda _{i}\xi _{i}^{2}$, with $\sum_{i=1}^{p}\xi
_{i}^{2}=\sum_{i=1}^{p}(\theta _{i}-\theta _{i}^{0})^{2}=\zeta ^{2}$
(see \citealp[p.~4605]{castilla2022}). Suppose that the negative eigenvalues $
\lambda _{i},i=1,...,p$ are numbered $\lambda _{p}\leq \lambda _{p-1}\leq
...\leq \lambda _{1}<0$. Then, $\sum_{i=1}^{p}\lambda _{i}\xi
_{i}^{2}\leq \lambda _{1}\sum_{i=1}^{p}\xi _{i}^{2}=\lambda
_{1}\zeta ^{2}$. Then, combining the first and second terms of~\eqref{c10}, we obtain
\begin{align*}
L_2
&=\frac12\sum_{j=1}^p\sum_{k=1}^p
\left\{-(a+1)I_{jk}(\boldsymbol{\theta}_0)(\theta_j-\theta_j^0)(\theta_k-\theta_k^0)\right\}\\
&\quad+\frac12\sum_{j=1}^p\sum_{k=1}^p
\left[\left(\frac{\partial^2 H_n^a(\boldsymbol{\theta})}{\partial\theta_j\partial\theta_k}\right)_{\boldsymbol{\theta}=\boldsymbol{\theta}_0}
+(a+1)I_{jk}(\boldsymbol{\theta}_0)\right]
(\theta_j-\theta_j^0)(\theta_k-\theta_k^0)\\
&\leq\sum_{i=1}^p\lambda_i\xi_i^2
+\left|\frac12\sum_{j=1}^p\sum_{k=1}^p
\left[\left(\frac{\partial^2 H_n^a(\boldsymbol{\theta})}{\partial\theta_j\partial\theta_k}\right)_{\boldsymbol{\theta}=\boldsymbol{\theta}_0}
+(a+1)I_{jk}(\boldsymbol{\theta}_0)\right]
(\theta_j-\theta_j^0)(\theta_k-\theta_k^0)\right|\\
&\leq\lambda_1\zeta^2+\frac12p^2\zeta^3,
\end{align*}
with probability tending to one, uniformly over the surface of $Q_\zeta$. The negative quadratic form is bounded above without taking its absolute value. Following the same sign argument as in \citet[p.~4605]{castilla2022}, since $\lambda_1<0$, the function $\lambda_1\zeta^2+\frac12p^2\zeta^3$ is negative for sufficiently small positive $\zeta$. More explicitly, take $c=-\lambda_1/2>0$ and $\zeta_0=-\lambda_1/p^2>0$. For $0<\zeta<\zeta_0$, we then have
\begin{equation}
L_2\leq-c\zeta^2,\qquad\text{with probability tending to }1.
\label{c11}
\end{equation}
This is the upper bound needed in the Taylor expansion of the objective.
To continue with $L_{3}$, condition (\textbf{C5}) and (\ref{c6}) give
\begin{align}
\left\vert L_{3}\right\vert  &=\left\vert \frac{1}{6}\sum\limits_{j=1}^{p}
\sum\limits_{k=1}^{p}\sum\limits_{\ell =1}^{p}\left( \frac{\partial
^{3}H_{n}^{a}(\boldsymbol{\theta })}{\partial \theta _{j}\partial \theta
_{k}\partial \theta _{\ell }}\right) _{\boldsymbol{\theta =\theta }^{\ast
}}(\theta _{j}-\theta _{j}^{0})(\theta _{k}-\theta _{k}^{0})(\theta _{\ell
}-\theta _{\ell }^{0})\right\vert   \notag \\
&\leq\frac{1}{6}\sum\limits_{j=1}^{p}\sum\limits_{k=1}^{p}\sum\limits_{
\ell =1}^{p}M_{jk\ell }\zeta ^{3}=b\zeta ^{3},  \label{c12}
\end{align}
with 
$
b=\frac{1}{6}\sum_{j=1}^{p}\sum_{k=1}^{p}\sum_{\ell
=1}^{p}M_{jk\ell }.
$
Based now on (\ref{c6}), (\ref{c8a}), (\ref{c11}) and (\ref{c12}) 
\begin{equation}
H_{n}^{a}(\boldsymbol{\theta })-H_{n}^{a}(\boldsymbol{\theta }
_{0})=L_{1}+L_{2}+L_{3}\leq p\zeta ^{3}-c\zeta ^{2}+b\zeta ^{3},  \label{c13}
\end{equation}
and the upper bound is less than zero if $\zeta <c/(b+p)$. Now, following \citet[p.~306]{basu2011} or \citet[p.~4605]{castilla2022}, the last inequality~\eqref{c13} implies that, for sufficiently small $\zeta$, there exists a sequence of interior maximizing roots
$
\widehat{\boldsymbol{\theta}}_{c,n}^a
=\widehat{\boldsymbol{\theta}}_{c,n}^a(\zeta)
$
such that
\begin{equation*}
\lim_{n\to\infty}P\left(\left\|\widehat{\boldsymbol{\theta}}_{c,n}^a(\zeta)
-\boldsymbol{\theta}_0\right\|_2<\zeta\right)=1,
\end{equation*}
where $\|\cdot\|_2$ represents the Euclidean norm. It remains to show that such a sequence can be selected independently of $\zeta$.

For this purpose, work in a fixed sufficiently small closed ball about $\boldsymbol{\theta}_0$, contained in the neighborhood in \textbf{(C2)--(C5)}. By~\eqref{c9a}, the Hessian at $\boldsymbol{\theta}_0$ converges to the negative definite matrix $-(a+1)I_a(\boldsymbol{\theta}_0)$. The uniform bound on its third derivatives in \textbf{(C5)} controls the change in every Hessian entry between $\boldsymbol{\theta}_0$ and another point of this ball. Choosing the ball sufficiently small, the Hessian is consequently negative definite throughout it with probability tending to one. On this event the criterion is strictly concave in the ball, so it has at most one stationary root there. The boundary argument already established gives an interior maximizing root in the ball with probability tending to one.

Equivalently, among the interior maximizing roots in this fixed ball, choose a root closest to $\boldsymbol{\theta}_0$. On the event just described it exists and is unique. A measurable arbitrary fallback can be used on the complement, whose probability tends to zero. For any smaller fixed $\zeta>0$, the boundary argument produces a maximizing root in $Q_\zeta$ with probability tending to one; uniqueness in the fixed ball makes it the same selected root. Therefore
\begin{equation*}
P\left(\left\|\widehat{\boldsymbol{\theta}}_{c,n}^a
-\boldsymbol{\theta}_0\right\|_2<\zeta\right)\longrightarrow1
\end{equation*}
for every $\zeta>0$. This concludes the proof of the existence of a sequence of consistent local maximizing solutions to the cumulative power divergence estimating equations, with probability tending to one. \hfill$\square$

\subsection{Proof of Theorem~\ref{thm:normality}}\label{proof:normality}

Following closely the exposition in \citet[p.~4605]{castilla2022}, \citet[p.~103]{leroy2016} and \citet[p.~306]{basu2011}, denote by
\begin{equation}
H_{n,j}^{a}(\boldsymbol{\theta })=\frac{\partial H_{n}^{a}(\boldsymbol{
\theta })}{\partial \theta _{j}},j=1,...,p,  \label{an1}
\end{equation}
with $\partial H_{n}^{a}(\boldsymbol{\theta })/\partial \theta _{j}$ defined
by (\ref{c7}). A Taylor series expansion of $H_{n,j}^{a}(\boldsymbol{\theta }
)$ around the true value $\boldsymbol{\theta }_{0}=(\theta
_{1}^{0},...,\theta _{p}^{0})^\top \in \Theta $ gives
\begin{align*}
H_{n,j}^{a}(\boldsymbol{\theta }) &=H_{n,j}^{a}(\boldsymbol{\theta }
_{0})+\sum\limits_{k=1}^{p}\left( \frac{\partial H_{n,j}^{a}(\boldsymbol{
\theta })}{\partial \theta _{k}}\right) _{\boldsymbol{\theta =\theta }
_{0}}(\theta _{k}-\theta _{k}^{0}) \\
&\quad+\frac{1}{2}\sum\limits_{k=1}^{p}\sum\limits_{\ell =1}^{p}\left( \frac{
\partial ^{2}H_{n,j}^{a}(\boldsymbol{\theta })}{\partial \theta _{k}\partial
\theta _{\ell }}\right) _{\boldsymbol{\theta =\theta }^{\ast }}(\theta
_{k}-\theta _{k}^{0})(\theta _{\ell }-\theta _{\ell }^{0}),
\end{align*}
or, taking into account (\ref{an1}), 
\begin{align}
H_{n,j}^{a}(\boldsymbol{\theta }) &=H_{n,j}^{a}(\boldsymbol{\theta }
_{0})+\sum\limits_{k=1}^{p}\left( \frac{\partial ^{2}H_{n}^{a}(\boldsymbol{
\theta })}{\partial \theta _{j}\partial \theta _{k}}\right) _{\boldsymbol{
\theta =\theta }_{0}}(\theta _{k}-\theta _{k}^{0})  \notag \\
&\quad+\frac{1}{2}\sum\limits_{k=1}^{p}\sum\limits_{\ell =1}^{p}\left( \frac{
\partial ^{3}H_{n}^{a}(\boldsymbol{\theta })}{\partial \theta _{j}\partial
\theta _{k}\partial \theta _{\ell }}\right) _{\boldsymbol{\theta =\theta }
^{\ast }}(\theta _{k}-\theta _{k}^{0})(\theta _{\ell }-\theta _{\ell }^{0}),
\label{an2}
\end{align}
with $\boldsymbol{\theta}^{\ast}$ on the segment connecting $\boldsymbol{\theta}$ and $\boldsymbol{\theta}_0$ (the intermediate point may depend on $j$), while assumptions \textbf{(C2)} and \textbf{(C3)} ensure the existence of the third order derivative in (\ref{an2}). It is
clear that the function $H_{n,j}^{a}(\boldsymbol{\theta })$ vanishes at the
MCumPDE $\widehat{\boldsymbol{\theta }}_{c,n}^{a}=\left( \widehat{\theta }
_{c,n,1}^{a},...,\widehat{\theta }_{c,n,p}^{a}\right) ^\top $, that is $
H_{n,j}^{a}(\widehat{\boldsymbol{\theta }}_{c,n}^{a})=0,$ and therefore
\begin{align}
H_{n,j}^{a}(\boldsymbol{\theta }_{0}) &=-\sum\limits_{k=1}^{p}\left( \frac{
\partial ^{2}H_{n}^{a}(\boldsymbol{\theta })}{\partial \theta _{j}\partial
\theta _{k}}\right) _{\boldsymbol{\theta =\theta }_{0}}(\widehat{\theta }
_{c,n,k}^{a}-\theta _{k}^{0})  \notag \\
&\quad-\frac{1}{2}\sum\limits_{k=1}^{p}\sum\limits_{\ell =1}^{p}\left( \frac{
\partial ^{3}H_{n}^{a}(\boldsymbol{\theta })}{\partial \theta _{j}\partial
\theta _{k}\partial \theta _{\ell }}\right) _{\boldsymbol{\theta =\theta }
^{\ast }}(\widehat{\theta }_{c,n,k}^{a}-\theta _{k}^{0})(\widehat{\theta }
_{c,n,\ell }^{a}-\theta _{\ell }^{0}).  \label{an3}
\end{align}
Based on this formulation, (\ref{an3}) is written as follows (see \citealp[p.~104]{leroy2016}),
\begin{equation}
\sqrt{n}\left( 
\begin{array}{c}
\frac{\partial H_{n}^{a}(\boldsymbol{\theta })}{\partial \theta _{1}} \\ 
\vdots  \\ 
\frac{\partial H_{n}^{a}(\boldsymbol{\theta })}{\partial \theta _{p}}
\end{array}
\right) _{\boldsymbol{\theta =\theta }_{0}}=\left( 
\begin{array}{ccc}
A_{11n} & \cdots  & A_{1pn} \\ 
\vdots  & \cdots  & \vdots  \\ 
A_{p1n} & \cdots  & A_{ppn}
\end{array}
\right) \left( 
\begin{array}{c}
\sqrt{n}\left( \widehat{\theta }_{c,n,1}^{a}-\theta _{1}^{0}\right)  \\ 
\vdots  \\ 
\sqrt{n}\left( \widehat{\theta }_{c,n,p}^{a}-\theta _{p}^{0}\right) 
\end{array}
\right) ,  \label{an4}
\end{equation}
with 
\begin{equation}
A_{jkn}=-\left( \frac{\partial ^{2}H_{n}^{a}(\boldsymbol{\theta })}{\partial
\theta _{j}\partial \theta _{k}}\right) _{\boldsymbol{\theta =\theta }_{0}}-
\frac{1}{2}\sum\limits_{\ell =1}^{p}\left( \frac{\partial ^{3}H_{n}^{a}(
\boldsymbol{\theta })}{\partial \theta _{j}\partial \theta _{k}\partial
\theta _{\ell }}\right) _{\boldsymbol{\theta =\theta }^{\ast }}\left( 
\widehat{\theta }_{c,n,\ell }^{a}-\theta _{\ell }^{0}\right) .  \label{an5}
\end{equation}
Relations (\ref{an4}) and (\ref{an5}) are expressed in matrix notation as
follows,
\begin{equation}
\boldsymbol{T}_{n}=\boldsymbol{A}_{n}\boldsymbol{Z}_{n},  \label{an6}
\end{equation}
with
$\boldsymbol{T}_{n}=\sqrt{n}\left( \frac{\partial H_{n}^{a}(\boldsymbol{
\theta })}{\partial \theta _{1}},...,\frac{\partial H_{n}^{a}(\boldsymbol{
\theta })}{\partial \theta _{p}}\right) _{\boldsymbol{\theta =\theta }
_{0}}^\top$, $\boldsymbol{A}_{n}=\left( A_{jkn}\right) _{j,k=1,...,p}$
and
$
\boldsymbol{Z}_{n}=\sqrt{n}\left( \widehat{\boldsymbol{\theta }}_{c,n}^{a}-
\boldsymbol{\theta }_{0}\right) =\sqrt{n}\left( \widehat{\theta }
_{c,n,1}^{a}-\theta _{1}^{0},...,\widehat{\theta }_{c,n,p}^{a}-\theta
_{p}^{0}\right) ^\top$.
Based now on the proof of consistency and on equation (\ref{c9a}), 
\begin{equation}
\left( \frac{\partial ^{2}H_{n}^{a}(\boldsymbol{\theta })}{\partial \theta
_{j}\partial \theta _{k}}\right) _{\boldsymbol{\theta =\theta }
_{0}}\rightarrow -(a+1)I_{jk}(\boldsymbol{\theta }_{0}),  \label{an7}
\end{equation}
with probability tending to $1$, as $n\rightarrow \infty $. Following closely the exposition in \citet[p.~4606]{castilla2022}, by assumption 
\textbf{(C5)}, $\left( \frac{\partial ^{3}H_{n}^{a}(\boldsymbol{\theta })}{
\partial \theta _{j}\partial \theta _{k}\partial \theta _{\ell }}\right) _{
\boldsymbol{\theta =\theta }^{\ast }}$ is bounded with probability tending
to one. Therefore, taking into account the consistency of $\widehat{
\boldsymbol{\theta }}_{c,n}^{a}$ we have that the second term of $A_{jkn}$
converges to zero in probability. Therefore, taking into account (\ref{an5})
and (\ref{an7}),
\begin{equation}
A_{jkn}\rightarrow (a+1)I_{jk}(\boldsymbol{\theta }_{0}),  \label{an8}
\end{equation}
with probability tending to one and $I_{jk}(\boldsymbol{\theta }_{0})$ is
defined by (\ref{c4}) of assumption \textbf{(C4)}.

Based on~\eqref{an6} and~\eqref{an8}, the random matrix $\boldsymbol{A}_n$ converges in probability to $(a+1)I_a(\boldsymbol{\theta}_0)$. Recall explicitly that $I_a(\boldsymbol{\theta}_0)$ is the $p\times p$ Fisher-type matrix with elements
\begin{equation*}
\begin{split}
I_a(\boldsymbol{\theta}_0)
&=\bigl(I_{jk}(\boldsymbol{\theta}_0)\bigr)_{j,k=1,\ldots,p}\\
&=\left(\int_{\D} F_0(\boldsymbol{x})^{1+a}
\left.\frac{\partial\log F_{\boldsymbol{\theta}}(\boldsymbol{x})}{\partial\theta_j}\right|_{\boldsymbol{\theta}=\boldsymbol{\theta}_0}
\left.\frac{\partial\log F_{\boldsymbol{\theta}}(\boldsymbol{x})}{\partial\theta_k}\right|_{\boldsymbol{\theta}=\boldsymbol{\theta}_0}
\dd\mu(\boldsymbol{x})\right)_{j,k=1,\ldots,p},\quad a>0.
\end{split}
\end{equation*}
Its positive definiteness in \textbf{(C4)} ensures that $\boldsymbol{A}_n$ is invertible with probability tending to one. On this event, the matrix relation~\eqref{an6} therefore gives the exact equation
\begin{equation}\label{an9}
\boldsymbol{Z}_n=\boldsymbol{A}_n^{-1}\boldsymbol{T}_n,
\qquad
\boldsymbol{A}_n^{-1}\overset{P}{\longrightarrow}(a+1)^{-1}I_a(\boldsymbol{\theta}_0)^{-1}.
\end{equation}
Thus the requested asymptotic distribution of $\boldsymbol{Z}_n$ will be the distributional limit of $(a+1)^{-1}I_a(\boldsymbol{\theta}_0)^{-1}\boldsymbol{T}_n$, once the latter limit has been established. This is an application of Slutsky's theorem, rather than a finite-sample equality in law. It remains to prove the asymptotic normality of the score vector
\begin{equation*}
\boldsymbol{T}_n=\sqrt n\left(
\frac{\partial H_n^a(\boldsymbol{\theta})}{\partial\theta_1},\ldots,
\frac{\partial H_n^a(\boldsymbol{\theta})}{\partial\theta_p}
\right)_{\boldsymbol{\theta}=\boldsymbol{\theta}_0}^{\top}.
\end{equation*}
In view of~\eqref{c7}, each coordinate may be written in the two equivalent forms
\begin{align}\label{an10}
\sqrt n\left.\frac{\partial H_n^a(\boldsymbol{\theta})}{\partial\theta_j}\right|_{\boldsymbol{\theta}=\boldsymbol{\theta}_0}
&=\int_{\D}(1+a)F_0(\boldsymbol{x})^{a-1}
\left.\frac{\partial F_{\boldsymbol{\theta}}(\boldsymbol{x})}{\partial\theta_j}\right|_{\boldsymbol{\theta}=\boldsymbol{\theta}_0}
\bigl[\sqrt n\{G_n(\boldsymbol{x})-F_0(\boldsymbol{x})\}\bigr]\dd\mu(\boldsymbol{x})\notag\\
&=\int_{\D}(1+a)F_0(\boldsymbol{x})^a
\left.\frac{\partial\log F_{\boldsymbol{\theta}}(\boldsymbol{x})}{\partial\theta_j}\right|_{\boldsymbol{\theta}=\boldsymbol{\theta}_0}
\bigl[\sqrt n\{G_n(\boldsymbol{x})-F_0(\boldsymbol{x})\}\bigr]\dd\mu(\boldsymbol{x}).
\end{align}
For the subsequent vector expressions, put
$
\boldsymbol{s}_F(\boldsymbol{x})
=\left.\nabla_{\boldsymbol{\theta}}\log F_{\boldsymbol{\theta}}(\boldsymbol{x})\right|_{\boldsymbol{\theta}=\boldsymbol{\theta}_0}.
$
Equation~\eqref{an10}, coordinate by coordinate, gives
\begin{equation}\label{an11}
\boldsymbol{T}_n=\int_{\D}(1+a)F_0(\boldsymbol{x})^a\boldsymbol{s}_F(\boldsymbol{x})
\bigl[\sqrt n\{G_n(\boldsymbol{x})-F_0(\boldsymbol{x})\}\bigr]\dd\mu(\boldsymbol{x}).
\end{equation}
Here $F_0(\boldsymbol{x})^a\boldsymbol{s}_F(\boldsymbol{x})=w_a(\boldsymbol{x};\boldsymbol{\theta}_0)$, and its norm is integrable by \textbf{(C3)}.

The distributional limit of the score is obtained from the multivariate empirical process. For a general multivariate distribution $F_0$, the empirical distribution process indexed by lower orthants satisfies
\begin{equation}\label{an:generalbridge}
\sqrt n\{G_n-F_0\}\Longrightarrow\mathbb{B}_{F_0}
\end{equation}
in the space of bounded functions on $\D$ with the supremum norm; see \citet{vandervaartwellner1996}. The limit is a centered Gaussian process with covariance kernel $K_0$ in~\eqref{eq:Kcdf}. No estimation of marginal distributions is involved in this limit. 

By \textbf{(C3)}, the linear map
$b\mapsto\int_D w_a(x;\theta_0)b(x)\,d\mu(x)$
is continuous for the supremum norm. Applying the continuous
mapping theorem to the empirical-process limit and the score
representation gives the following joint Gaussian limit
%
\begin{equation}\label{an15}
\boldsymbol{T}_n\overset{\mathcal{L}}{\longrightarrow}
\int_{\D}\mathbb{B}_{F_0}(\boldsymbol{x})
\bigl[(1+a)F_0(\boldsymbol{x})^a\boldsymbol{s}_F(\boldsymbol{x})\bigr]\dd\mu(\boldsymbol{x}).
\end{equation}
For continuous margins this is equivalently the vector integral with $\mathbb{B}_{C_0}\{\boldsymbol{F}(\boldsymbol{x})\}$ displayed in the coordinatewise calculation. For arbitrary margins, including discontinuous ones, \eqref{an:generalbridge} directly yields~\eqref{an15}; no copula representation is required.

The vector on the right-hand side of~\eqref{an15} is Gaussian, with mean $\boldsymbol{0}_p$. Its variance-covariance matrix is obtained by integrating the covariance kernel of the limiting process against the two score weights:
\begin{align}\label{an16}
\boldsymbol{\Sigma}_T
&=(1+a)^2\int_{\D}\int_{\D}
\bigl[F_0(\boldsymbol{x}\wedge\boldsymbol{y})-F_0(\boldsymbol{x})F_0(\boldsymbol{y})\bigr]
\bigl[F_0(\boldsymbol{x})^a\boldsymbol{s}_F(\boldsymbol{x})\bigr]\notag\\
&\quad\times\bigl[F_0(\boldsymbol{y})^a\boldsymbol{s}_F(\boldsymbol{y})\bigr]^\top 
\dd\mu(\boldsymbol{x})\dd\mu(\boldsymbol{y})
=(1+a)^2 B_a(\boldsymbol{\theta}_0).
\end{align}
In the copula representation of the same covariance, this becomes
\begin{align*}
\boldsymbol{\Sigma}_T
&=(1+a)^2\int_{\D}\int_{\D}
\bigl[C_0\{\boldsymbol{F}(\boldsymbol{x})\wedge\boldsymbol{F}(\boldsymbol{y})\}
-C_0\{\boldsymbol{F}(\boldsymbol{x})\}C_0\{\boldsymbol{F}(\boldsymbol{y})\}\bigr]\\
&\qquad\qquad\times\bigl[F_0(\boldsymbol{x})^a\boldsymbol{s}_F(\boldsymbol{x})\bigr]
\bigl[F_0(\boldsymbol{y})^a\boldsymbol{s}_F(\boldsymbol{y})\bigr]^\top 
\dd\mu(\boldsymbol{x})\dd\mu(\boldsymbol{y}),
\end{align*}
where $\boldsymbol{F}(\boldsymbol{y})=(F_1(y_1),\ldots,F_m(y_m))^\top $. The products of the column vectors are outer products. The equality of these two covariance representations follows from
$
\boldsymbol{F}(\boldsymbol{x}\wedge\boldsymbol{y})
=\boldsymbol{F}(\boldsymbol{x})\wedge\boldsymbol{F}(\boldsymbol{y})$,
$F_0=C_0\circ\boldsymbol{F}$.
The integrable score weight in \textbf{(C3)} and the bounded covariance kernel justify the double integral. The factor $(1+a)^2$ is part of the score covariance; it will cancel in the covariance of the estimator.
Then, in view of~\eqref{an9} and~\eqref{an15}, Slutsky's theorem gives
\begin{align}\label{an17}
\boldsymbol{Z}_n
=\sqrt n\bigl(\widehat{\boldsymbol{\theta}}_{c,n}^a-\boldsymbol{\theta}_0\bigr)
=\boldsymbol{A}_n^{-1}\boldsymbol{T}_n
\overset{\mathcal{L}}{\longrightarrow}
I_a(\boldsymbol{\theta}_0)^{-1}\int_{\D}
\mathbb{B}_{F_0}(\boldsymbol{x})\bigl[F_0(\boldsymbol{x})^a\boldsymbol{s}_F(\boldsymbol{x})\bigr]\dd\mu(\boldsymbol{x}).
\end{align}
For continuous margins the right-hand side has the equivalent copula representation
$
I_a(\boldsymbol{\theta}_0)^{-1}\int_{\D}
\mathbb{B}_{C_0}\{\boldsymbol{F}(\boldsymbol{x})\}
\bigl[F_0(\boldsymbol{x})^a\boldsymbol{s}_F(\boldsymbol{x})\bigr]\dd\mu(\boldsymbol{x})$.
Based on this analysis, and specifically on~\eqref{an15}, \eqref{an16} and~\eqref{an17}, we conclude that
\begin{equation*}
\sqrt n\bigl(\widehat{\boldsymbol{\theta}}_{c,n}^a-\boldsymbol{\theta}_0\bigr)
\overset{\mathcal{L}}{\longrightarrow}N_p(\boldsymbol{0}_p,\boldsymbol{\Sigma}_a),
\end{equation*}
with
$
\boldsymbol{\Sigma}_a
=(1+a)^{-2}I_a(\boldsymbol{\theta}_0)^{-1}
\boldsymbol{\Sigma}_T I_a(\boldsymbol{\theta}_0)^{-1}
=I_a(\boldsymbol{\theta}_0)^{-1}B_a(\boldsymbol{\theta}_0)I_a(\boldsymbol{\theta}_0)^{-1}.
$
This proves Theorem~\ref{thm:normality}. \hfill$\square$

\end{appendices}


\bibliography{sn-bibliography}

\end{document}


\title[Supplementary Material to: Minimum cumulative power divergence estimation with applications to copulas]{Supplementary Material to: Minimum cumulative power divergence estimation with applications to copulas}


%
%
%

\author*[1]{\fnm{Ostap} \sur{Okhrin}}
\email{ostap.okhrin@tu-dresden.de}

\author[2]{\fnm{Konstantinos} \sur{Zografos}}
\email{kzograf@uoi.gr}

\affil*[1]{%
  \orgdiv{Chair of Applied Statistics,
    Faculty of Transport and Traffic Sciences},
  \orgname{Technische Universit\"at Dresden},
  \orgaddress{%
    \city{Dresden},
    \country{Germany}}}

\affil[2]{%
  \orgdiv{Department of Mathematics},
  \orgname{University of Ioannina},
  \orgaddress{%
    \city{Ioannina},
    \country{Greece}}}




\maketitle

\section*{Scope and interpretation}
This supplement provides the complete numerical comparisons for the simulation
study in Section~5 of the paper. The design comprises 216 cells: two sample
sizes ($n=100,1000$), four uncontaminated dependence levels
($\tau_0=0.1,0.3,0.5,0.9$), nine generating/fitted-family combinations
(Gumbel, Clayton and Frank), and three independence-contamination fractions
($\varepsilon=0,0.1,0.5$). Each reported entry is based on 200 replications.

Tables~\ref{tab:comparison-n1000-low}--\ref{tab:comparison-n1000-high} (\ref{tab:comparison-n100-low}--\ref{tab:comparison-n100-high})
report the $n=1000$ ($n=100$) comparison of MCopPDE and EK--$\beta$ at the common power
$a=q=0.5$.
Tables~\ref{tab:complete-n100-gumbel-gumbel}--\ref{tab:complete-n1000-frank-frank}
report mean absolute error (MAE), signed bias and the Monte Carlo standard
error (MCSE) of MAE for all 216 cells and all evaluated powers of MCopPDE,
EK--$\alpha$, EK--$\beta$ and EK--$\gamma$, together with sample Kendall tau.
The latter is independent of the fitted family and is repeated as a reference.
The tables focus on these estimators; likelihood benchmarks are discussed
separately in Section~5.3 of the paper.

All errors refer to the uncontaminated Kendall tau $\tau_0$. At positive
contamination they assess recovery of the clean component, rather than
estimation of the Kendall tau of the contaminated mixture. For replication
$b$, write $E_b=|\widehat\tau_b-\tau_0|$. Then
\[
\operatorname{MAE}=\frac{1}{200}\sum_{b=1}^{200}E_b,\qquad
\operatorname{Bias}=\frac{1}{200}\sum_{b=1}^{200}(\widehat\tau_b-\tau_0),\qquad
\operatorname{MCSE}(\operatorname{MAE})=\frac{s_E}{\sqrt{200}},
\]
where $s_E$ is the sample standard deviation of the absolute errors.
The columns labelled SE display this MCSE in units of $10^{-4}$.

\clearpage
\section{\texorpdfstring{Systematic comparison for $n=1000$}{Systematic comparison for n=1000}}

\begin{table}[!htbp]
\centering\small
\caption{Systematic comparison at $n=1000$ and $a=q=0.5$. Each entry gives mean absolute error above, with signed bias in square brackets below; both refer to the uncontaminated $\tau_0$. All nine generating/fitted-family combinations and all three contamination fractions are included for each displayed dependence level. Bold family pairs identify correct specification. Every entry uses 200 replications. MCSEs and results for all other evaluated powers and divergence estimators are reported in Tables~\ref{tab:complete-n1000-gumbel-gumbel}--\ref{tab:complete-n1000-frank-frank}.}
\label{tab:comparison-n1000-low}
\setlength{\tabcolsep}{3.5pt}
\renewcommand{\arraystretch}{1.05}
\begin{tabular}{@{}ll*{6}{r}@{}}
\toprule
& & \multicolumn{2}{c}{$\varepsilon=0$} & \multicolumn{2}{c}{$\varepsilon=0.1$} & \multicolumn{2}{c}{$\varepsilon=0.5$}\\
\cmidrule(lr){3-4}\cmidrule(lr){5-6}\cmidrule(lr){7-8}
Generating & Fitted & MCopPDE & EK--$\beta$ & MCopPDE & EK--$\beta$ & MCopPDE & EK--$\beta$\\
\midrule
\multicolumn{8}{@{}l}{\textbf{Kendall tau $\tau_0=0.1$}}\\[3pt]
\textbf{Gumbel} & \textbf{Gumbel} & \shortstack{$0.0174$\\$[+0.0017]$} & \shortstack{$0.0182$\\$[+0.0063]$} & \shortstack{$0.0206$\\$[-0.0084]$} & \shortstack{$0.0195$\\$[-0.0037]$} & \shortstack{$0.0513$\\$[-0.0509]$} & \shortstack{$0.0468$\\$[-0.0462]$}\\[2pt]
Gumbel & Clayton & \shortstack{$0.0249$\\$[-0.0222]$} & \shortstack{$0.0227$\\$[-0.0194]$} & \shortstack{$0.0317$\\$[-0.0302]$} & \shortstack{$0.0292$\\$[-0.0272]$} & \shortstack{$0.0634$\\$[-0.0634]$} & \shortstack{$0.0599$\\$[-0.0599]$}\\[2pt]
Gumbel & Frank & \shortstack{$0.0175$\\$[-0.0077]$} & \shortstack{$0.0167$\\$[-0.0037]$} & \shortstack{$0.0234$\\$[-0.0169]$} & \shortstack{$0.0215$\\$[-0.0129]$} & \shortstack{$0.0559$\\$[-0.0557]$} & \shortstack{$0.0520$\\$[-0.0516]$}\\[2pt]
\addlinespace[4pt]
Clayton & Gumbel & \shortstack{$0.0204$\\$[+0.0103]$} & \shortstack{$0.0223$\\$[+0.0146]$} & \shortstack{$0.0187$\\$[-0.0025]$} & \shortstack{$0.0185$\\$[+0.0020]$} & \shortstack{$0.0478$\\$[-0.0477]$} & \shortstack{$0.0432$\\$[-0.0430]$}\\[2pt]
\textbf{Clayton} & \textbf{Clayton} & \shortstack{$0.0156$\\$[+0.0004]$} & \shortstack{$0.0158$\\$[+0.0043]$} & \shortstack{$0.0176$\\$[-0.0111]$} & \shortstack{$0.0160$\\$[-0.0072]$} & \shortstack{$0.0521$\\$[-0.0521]$} & \shortstack{$0.0482$\\$[-0.0482]$}\\[2pt]
Clayton & Frank & \shortstack{$0.0176$\\$[+0.0052]$} & \shortstack{$0.0189$\\$[+0.0094]$} & \shortstack{$0.0179$\\$[-0.0068]$} & \shortstack{$0.0174$\\$[-0.0026]$} & \shortstack{$0.0498$\\$[-0.0497]$} & \shortstack{$0.0456$\\$[-0.0455]$}\\[2pt]
\addlinespace[4pt]
Frank & Gumbel & \shortstack{$0.0174$\\$[+0.0061]$} & \shortstack{$0.0185$\\$[+0.0106]$} & \shortstack{$0.0176$\\$[-0.0054]$} & \shortstack{$0.0171$\\$[-0.0008]$} & \shortstack{$0.0477$\\$[-0.0474]$} & \shortstack{$0.0431$\\$[-0.0426]$}\\[2pt]
Frank & Clayton & \shortstack{$0.0189$\\$[-0.0143]$} & \shortstack{$0.0173$\\$[-0.0119]$} & \shortstack{$0.0257$\\$[-0.0237]$} & \shortstack{$0.0234$\\$[-0.0209]$} & \shortstack{$0.0578$\\$[-0.0578]$} & \shortstack{$0.0543$\\$[-0.0542]$}\\[2pt]
\textbf{Frank} & \textbf{Frank} & \shortstack{$0.0158$\\$[-0.0012]$} & \shortstack{$0.0158$\\$[+0.0031]$} & \shortstack{$0.0187$\\$[-0.0118]$} & \shortstack{$0.0171$\\$[-0.0076]$} & \shortstack{$0.0511$\\$[-0.0510]$} & \shortstack{$0.0469$\\$[-0.0467]$}\\[2pt]
\addlinespace[6pt]
\multicolumn{8}{@{}l}{\textbf{Kendall tau $\tau_0=0.3$}}\\[3pt]
\textbf{Gumbel} & \textbf{Gumbel} & \shortstack{$0.0169$\\$[+0.0000]$} & \shortstack{$0.0173$\\$[+0.0044]$} & \shortstack{$0.0321$\\$[-0.0312]$} & \shortstack{$0.0280$\\$[-0.0265]$} & \shortstack{$0.1493$\\$[-0.1493]$} & \shortstack{$0.1441$\\$[-0.1441]$}\\[2pt]
Gumbel & Clayton & \shortstack{$0.0489$\\$[-0.0489]$} & \shortstack{$0.0516$\\$[-0.0516]$} & \shortstack{$0.0781$\\$[-0.0781]$} & \shortstack{$0.0795$\\$[-0.0795]$} & \shortstack{$0.1802$\\$[-0.1802]$} & \shortstack{$0.1778$\\$[-0.1778]$}\\[2pt]
Gumbel & Frank & \shortstack{$0.0251$\\$[-0.0229]$} & \shortstack{$0.0222$\\$[-0.0183]$} & \shortstack{$0.0523$\\$[-0.0523]$} & \shortstack{$0.0477$\\$[-0.0477]$} & \shortstack{$0.1620$\\$[-0.1620]$} & \shortstack{$0.1573$\\$[-0.1573]$}\\[2pt]
\addlinespace[4pt]
Clayton & Gumbel & \shortstack{$0.0372$\\$[+0.0365]$} & \shortstack{$0.0464$\\$[+0.0462]$} & \shortstack{$0.0160$\\$[-0.0005]$} & \shortstack{$0.0179$\\$[+0.0085]$} & \shortstack{$0.1377$\\$[-0.1377]$} & \shortstack{$0.1316$\\$[-0.1316]$}\\[2pt]
\textbf{Clayton} & \textbf{Clayton} & \shortstack{$0.0142$\\$[+0.0021]$} & \shortstack{$0.0148$\\$[+0.0057]$} & \shortstack{$0.0316$\\$[-0.0311]$} & \shortstack{$0.0277$\\$[-0.0270]$} & \shortstack{$0.1543$\\$[-0.1543]$} & \shortstack{$0.1495$\\$[-0.1495]$}\\[2pt]
Clayton & Frank & \shortstack{$0.0181$\\$[+0.0127]$} & \shortstack{$0.0250$\\$[+0.0226]$} & \shortstack{$0.0229$\\$[-0.0208]$} & \shortstack{$0.0179$\\$[-0.0118]$} & \shortstack{$0.1467$\\$[-0.1467]$} & \shortstack{$0.1408$\\$[-0.1408]$}\\[2pt]
\addlinespace[4pt]
Frank & Gumbel & \shortstack{$0.0230$\\$[+0.0189]$} & \shortstack{$0.0243$\\$[+0.0215]$} & \shortstack{$0.0210$\\$[-0.0129]$} & \shortstack{$0.0194$\\$[-0.0095]$} & \shortstack{$0.1428$\\$[-0.1428]$} & \shortstack{$0.1379$\\$[-0.1379]$}\\[2pt]
Frank & Clayton & \shortstack{$0.0275$\\$[-0.0244]$} & \shortstack{$0.0316$\\$[-0.0299]$} & \shortstack{$0.0546$\\$[-0.0546]$} & \shortstack{$0.0582$\\$[-0.0582]$} & \shortstack{$0.1707$\\$[-0.1707]$} & \shortstack{$0.1689$\\$[-0.1689]$}\\[2pt]
\textbf{Frank} & \textbf{Frank} & \shortstack{$0.0172$\\$[-0.0031]$} & \shortstack{$0.0167$\\$[+0.0010]$} & \shortstack{$0.0342$\\$[-0.0328]$} & \shortstack{$0.0303$\\$[-0.0283]$} & \shortstack{$0.1537$\\$[-0.1537]$} & \shortstack{$0.1489$\\$[-0.1489]$}\\[2pt]
\bottomrule
\end{tabular}
\end{table}

\begin{table}[!htbp]
\centering\small
\caption{Systematic comparison at $n=1000$ and $a=q=0.5$. Each entry gives mean absolute error above, with signed bias in square brackets below; both refer to the uncontaminated $\tau_0$. All nine generating/fitted-family combinations and all three contamination fractions are included for each displayed dependence level. Bold family pairs identify correct specification. Every entry uses 200 replications. MCSEs and results for all other evaluated powers and divergence estimators are reported in Tables~\ref{tab:complete-n1000-gumbel-gumbel}--\ref{tab:complete-n1000-frank-frank}.}
\label{tab:comparison-n1000-high}
\setlength{\tabcolsep}{3.5pt}
\renewcommand{\arraystretch}{1.05}
\begin{tabular}{@{}ll*{6}{r}@{}}
\toprule
& & \multicolumn{2}{c}{$\varepsilon=0$} & \multicolumn{2}{c}{$\varepsilon=0.1$} & \multicolumn{2}{c}{$\varepsilon=0.5$}\\
\cmidrule(lr){3-4}\cmidrule(lr){5-6}\cmidrule(lr){7-8}
Generating & Fitted & MCopPDE & EK--$\beta$ & MCopPDE & EK--$\beta$ & MCopPDE & EK--$\beta$\\
\midrule
\multicolumn{8}{@{}l}{\textbf{Kendall tau $\tau_0=0.5$}}\\[3pt]
\textbf{Gumbel} & \textbf{Gumbel} & \shortstack{$0.0148$\\$[-0.0018]$} & \shortstack{$0.0145$\\$[+0.0024]$} & \shortstack{$0.0547$\\$[-0.0547]$} & \shortstack{$0.0485$\\$[-0.0485]$} & \shortstack{$0.2604$\\$[-0.2604]$} & \shortstack{$0.2524$\\$[-0.2524]$}\\[2pt]
Gumbel & Clayton & \shortstack{$0.0401$\\$[-0.0399]$} & \shortstack{$0.0427$\\$[-0.0426]$} & \shortstack{$0.0996$\\$[-0.0996]$} & \shortstack{$0.1010$\\$[-0.1010]$} & \shortstack{$0.3017$\\$[-0.3017]$} & \shortstack{$0.2985$\\$[-0.2985]$}\\[2pt]
Gumbel & Frank & \shortstack{$0.0324$\\$[-0.0322]$} & \shortstack{$0.0252$\\$[-0.0242]$} & \shortstack{$0.0842$\\$[-0.0842]$} & \shortstack{$0.0753$\\$[-0.0753]$} & \shortstack{$0.2790$\\$[-0.2790]$} & \shortstack{$0.2715$\\$[-0.2715]$}\\[2pt]
\addlinespace[4pt]
Clayton & Gumbel & \shortstack{$0.0522$\\$[+0.0522]$} & \shortstack{$0.0807$\\$[+0.0807]$} & \shortstack{$0.0157$\\$[-0.0100]$} & \shortstack{$0.0197$\\$[+0.0165]$} & \shortstack{$0.2395$\\$[-0.2395]$} & \shortstack{$0.2253$\\$[-0.2253]$}\\[2pt]
\textbf{Clayton} & \textbf{Clayton} & \shortstack{$0.0126$\\$[-0.0015]$} & \shortstack{$0.0116$\\$[+0.0025]$} & \shortstack{$0.0576$\\$[-0.0576]$} & \shortstack{$0.0522$\\$[-0.0522]$} & \shortstack{$0.2692$\\$[-0.2692]$} & \shortstack{$0.2608$\\$[-0.2608]$}\\[2pt]
Clayton & Frank & \shortstack{$0.0146$\\$[+0.0106]$} & \shortstack{$0.0380$\\$[+0.0379]$} & \shortstack{$0.0468$\\$[-0.0468]$} & \shortstack{$0.0240$\\$[-0.0225]$} & \shortstack{$0.2567$\\$[-0.2567]$} & \shortstack{$0.2438$\\$[-0.2438]$}\\[2pt]
\addlinespace[4pt]
Frank & Gumbel & \shortstack{$0.0328$\\$[+0.0326]$} & \shortstack{$0.0327$\\$[+0.0326]$} & \shortstack{$0.0263$\\$[-0.0245]$} & \shortstack{$0.0233$\\$[-0.0209]$} & \shortstack{$0.2416$\\$[-0.2416]$} & \shortstack{$0.2338$\\$[-0.2338]$}\\[2pt]
Frank & Clayton & \shortstack{$0.0151$\\$[-0.0013]$} & \shortstack{$0.0166$\\$[-0.0095]$} & \shortstack{$0.0647$\\$[-0.0647]$} & \shortstack{$0.0701$\\$[-0.0701]$} & \shortstack{$0.2814$\\$[-0.2814]$} & \shortstack{$0.2793$\\$[-0.2793]$}\\[2pt]
\textbf{Frank} & \textbf{Frank} & \shortstack{$0.0130$\\$[-0.0003]$} & \shortstack{$0.0130$\\$[+0.0038]$} & \shortstack{$0.0554$\\$[-0.0554]$} & \shortstack{$0.0488$\\$[-0.0488]$} & \shortstack{$0.2599$\\$[-0.2599]$} & \shortstack{$0.2520$\\$[-0.2520]$}\\[2pt]
\addlinespace[6pt]
\multicolumn{8}{@{}l}{\textbf{Kendall tau $\tau_0=0.9$}}\\[3pt]
\textbf{Gumbel} & \textbf{Gumbel} & \shortstack{$0.0035$\\$[-0.0011]$} & \shortstack{$0.0040$\\$[+0.0030]$} & \shortstack{$0.1243$\\$[-0.1243]$} & \shortstack{$0.0900$\\$[-0.0900]$} & \shortstack{$0.5214$\\$[-0.5214]$} & \shortstack{$0.4639$\\$[-0.4639]$}\\[2pt]
Gumbel & Clayton & \shortstack{$0.0204$\\$[+0.0204]$} & \shortstack{$0.0392$\\$[+0.0392]$} & \shortstack{$0.1199$\\$[-0.1199]$} & \shortstack{$0.0790$\\$[-0.0790]$} & \shortstack{$0.5668$\\$[-0.5668]$} & \shortstack{$0.5196$\\$[-0.5196]$}\\[2pt]
Gumbel & Frank & \shortstack{$0.0091$\\$[-0.0091]$} & \shortstack{$0.0070$\\$[+0.0068]$} & \shortstack{$0.1481$\\$[-0.1481]$} & \shortstack{$0.0996$\\$[-0.0996]$} & \shortstack{$0.5489$\\$[-0.5489]$} & \shortstack{$0.4937$\\$[-0.4937]$}\\[2pt]
\addlinespace[4pt]
Clayton & Gumbel & \shortstack{$0.0297$\\$[+0.0297]$} & \shortstack{$0.0473$\\$[+0.0473]$} & \shortstack{$0.1061$\\$[-0.1061]$} & \shortstack{$0.0498$\\$[-0.0498]$} & \shortstack{$0.5123$\\$[-0.5123]$} & \shortstack{$0.4428$\\$[-0.4428]$}\\[2pt]
\textbf{Clayton} & \textbf{Clayton} & \shortstack{$0.0042$\\$[-0.0011]$} & \shortstack{$0.0052$\\$[+0.0039]$} & \shortstack{$0.1191$\\$[-0.1191]$} & \shortstack{$0.0941$\\$[-0.0941]$} & \shortstack{$0.5566$\\$[-0.5566]$} & \shortstack{$0.4996$\\$[-0.4996]$}\\[2pt]
Clayton & Frank & \shortstack{$0.0100$\\$[+0.0100]$} & \shortstack{$0.0463$\\$[+0.0463]$} & \shortstack{$0.1334$\\$[-0.1334]$} & \shortstack{$0.0660$\\$[-0.0660]$} & \shortstack{$0.5403$\\$[-0.5403]$} & \shortstack{$0.4751$\\$[-0.4751]$}\\[2pt]
\addlinespace[4pt]
Frank & Gumbel & \shortstack{$0.0148$\\$[+0.0148]$} & \shortstack{$0.0171$\\$[+0.0171]$} & \shortstack{$0.1147$\\$[-0.1147]$} & \shortstack{$0.0775$\\$[-0.0775]$} & \shortstack{$0.5175$\\$[-0.5175]$} & \shortstack{$0.4558$\\$[-0.4558]$}\\[2pt]
Frank & Clayton & \shortstack{$0.0219$\\$[+0.0219]$} & \shortstack{$0.0252$\\$[+0.0252]$} & \shortstack{$0.1121$\\$[-0.1121]$} & \shortstack{$0.0794$\\$[-0.0794]$} & \shortstack{$0.5626$\\$[-0.5626]$} & \shortstack{$0.5104$\\$[-0.5104]$}\\[2pt]
\textbf{Frank} & \textbf{Frank} & \shortstack{$0.0028$\\$[-0.0013]$} & \shortstack{$0.0043$\\$[+0.0037]$} & \shortstack{$0.1396$\\$[-0.1396]$} & \shortstack{$0.0955$\\$[-0.0955]$} & \shortstack{$0.5451$\\$[-0.5451]$} & \shortstack{$0.4857$\\$[-0.4857]$}\\[2pt]
\bottomrule
\end{tabular}
\end{table}
\clearpage
\section{\texorpdfstring{Systematic comparison for $n=100$}{Systematic comparison for n=100}}

\begin{table}[!htbp]
\centering\small
\caption{Systematic comparison at $n=100$ and $a=q=0.5$. Each entry gives mean absolute error above, with signed bias in square brackets below; both refer to the uncontaminated $\tau_0$. All nine generating/fitted-family combinations and all three contamination fractions are included for each displayed dependence level. Bold family pairs identify correct specification. Every entry uses 200 replications. MCSEs and results for all other evaluated powers and divergence estimators are reported in Tables~\ref{tab:complete-n100-gumbel-gumbel}--\ref{tab:complete-n100-frank-frank}.}
\label{tab:comparison-n100-low}
\setlength{\tabcolsep}{3.5pt}
\renewcommand{\arraystretch}{1.05}
\begin{tabular}{@{}ll*{6}{r}@{}}
\toprule
& & \multicolumn{2}{c}{$\varepsilon=0$} & \multicolumn{2}{c}{$\varepsilon=0.1$} & \multicolumn{2}{c}{$\varepsilon=0.5$}\\
\cmidrule(lr){3-4}\cmidrule(lr){5-6}\cmidrule(lr){7-8}
Generating & Fitted & MCopPDE & EK--$\beta$ & MCopPDE & EK--$\beta$ & MCopPDE & EK--$\beta$\\
\midrule
\multicolumn{8}{@{}l}{\textbf{Kendall tau $\tau_0=0.1$}}\\[3pt]
\textbf{Gumbel} & \textbf{Gumbel} & \shortstack{$0.0546$\\$[+0.0027]$} & \shortstack{$0.0711$\\$[+0.0521]$} & \shortstack{$0.0560$\\$[-0.0060]$} & \shortstack{$0.0675$\\$[+0.0418]$} & \shortstack{$0.0566$\\$[-0.0344]$} & \shortstack{$0.0523$\\$[+0.0108]$}\\[2pt]
Gumbel & Clayton & \shortstack{$0.0513$\\$[-0.0234]$} & \shortstack{$0.0506$\\$[+0.0193]$} & \shortstack{$0.0538$\\$[-0.0304]$} & \shortstack{$0.0503$\\$[+0.0118]$} & \shortstack{$0.0618$\\$[-0.0514]$} & \shortstack{$0.0436$\\$[-0.0119]$}\\[2pt]
Gumbel & Frank & \shortstack{$0.0518$\\$[-0.0077]$} & \shortstack{$0.0614$\\$[+0.0371]$} & \shortstack{$0.0539$\\$[-0.0155]$} & \shortstack{$0.0596$\\$[+0.0280]$} & \shortstack{$0.0578$\\$[-0.0411]$} & \shortstack{$0.0485$\\$[-0.0001]$}\\[2pt]
\addlinespace[4pt]
Clayton & Gumbel & \shortstack{$0.0588$\\$[+0.0098]$} & \shortstack{$0.0782$\\$[+0.0581]$} & \shortstack{$0.0573$\\$[-0.0007]$} & \shortstack{$0.0700$\\$[+0.0482]$} & \shortstack{$0.0605$\\$[-0.0364]$} & \shortstack{$0.0559$\\$[+0.0082]$}\\[2pt]
\textbf{Clayton} & \textbf{Clayton} & \shortstack{$0.0515$\\$[-0.0054]$} & \shortstack{$0.0577$\\$[+0.0378]$} & \shortstack{$0.0499$\\$[-0.0146]$} & \shortstack{$0.0521$\\$[+0.0290]$} & \shortstack{$0.0587$\\$[-0.0466]$} & \shortstack{$0.0482$\\$[-0.0070]$}\\[2pt]
Clayton & Frank & \shortstack{$0.0548$\\$[+0.0025]$} & \shortstack{$0.0687$\\$[+0.0466]$} & \shortstack{$0.0538$\\$[-0.0072]$} & \shortstack{$0.0619$\\$[+0.0373]$} & \shortstack{$0.0595$\\$[-0.0408]$} & \shortstack{$0.0524$\\$[-0.0006]$}\\[2pt]
\addlinespace[4pt]
Frank & Gumbel & \shortstack{$0.0551$\\$[+0.0057]$} & \shortstack{$0.0713$\\$[+0.0534]$} & \shortstack{$0.0569$\\$[-0.0044]$} & \shortstack{$0.0682$\\$[+0.0420]$} & \shortstack{$0.0637$\\$[-0.0448]$} & \shortstack{$0.0554$\\$[-0.0050]$}\\[2pt]
Frank & Clayton & \shortstack{$0.0501$\\$[-0.0180]$} & \shortstack{$0.0507$\\$[+0.0239]$} & \shortstack{$0.0527$\\$[-0.0257]$} & \shortstack{$0.0518$\\$[+0.0148]$} & \shortstack{$0.0670$\\$[-0.0589]$} & \shortstack{$0.0509$\\$[-0.0241]$}\\[2pt]
\textbf{Frank} & \textbf{Frank} & \shortstack{$0.0519$\\$[-0.0029]$} & \shortstack{$0.0628$\\$[+0.0409]$} & \shortstack{$0.0540$\\$[-0.0121]$} & \shortstack{$0.0614$\\$[+0.0304]$} & \shortstack{$0.0644$\\$[-0.0499]$} & \shortstack{$0.0537$\\$[-0.0141]$}\\[2pt]
\addlinespace[6pt]
\multicolumn{8}{@{}l}{\textbf{Kendall tau $\tau_0=0.3$}}\\[3pt]
\textbf{Gumbel} & \textbf{Gumbel} & \shortstack{$0.0566$\\$[-0.0040]$} & \shortstack{$0.0693$\\$[+0.0468]$} & \shortstack{$0.0657$\\$[-0.0335]$} & \shortstack{$0.0611$\\$[+0.0178]$} & \shortstack{$0.1513$\\$[-0.1512]$} & \shortstack{$0.1027$\\$[-0.0987]$}\\[2pt]
Gumbel & Clayton & \shortstack{$0.0726$\\$[-0.0526]$} & \shortstack{$0.0567$\\$[-0.0091]$} & \shortstack{$0.0901$\\$[-0.0804]$} & \shortstack{$0.0641$\\$[-0.0365]$} & \shortstack{$0.1848$\\$[-0.1848]$} & \shortstack{$0.1389$\\$[-0.1389]$}\\[2pt]
Gumbel & Frank & \shortstack{$0.0596$\\$[-0.0270]$} & \shortstack{$0.0579$\\$[+0.0212]$} & \shortstack{$0.0723$\\$[-0.0548]$} & \shortstack{$0.0575$\\$[-0.0066]$} & \shortstack{$0.1651$\\$[-0.1651]$} & \shortstack{$0.1179$\\$[-0.1170]$}\\[2pt]
\addlinespace[4pt]
Clayton & Gumbel & \shortstack{$0.0619$\\$[+0.0287]$} & \shortstack{$0.0963$\\$[+0.0872]$} & \shortstack{$0.0572$\\$[-0.0018]$} & \shortstack{$0.0738$\\$[+0.0556]$} & \shortstack{$0.1442$\\$[-0.1416]$} & \shortstack{$0.0993$\\$[-0.0899]$}\\[2pt]
\textbf{Clayton} & \textbf{Clayton} & \shortstack{$0.0531$\\$[-0.0051]$} & \shortstack{$0.0603$\\$[+0.0400]$} & \shortstack{$0.0593$\\$[-0.0334]$} & \shortstack{$0.0519$\\$[+0.0121]$} & \shortstack{$0.1611$\\$[-0.1604]$} & \shortstack{$0.1180$\\$[-0.1155]$}\\[2pt]
Clayton & Frank & \shortstack{$0.0537$\\$[+0.0056]$} & \shortstack{$0.0740$\\$[+0.0591]$} & \shortstack{$0.0559$\\$[-0.0225]$} & \shortstack{$0.0591$\\$[+0.0298]$} & \shortstack{$0.1529$\\$[-0.1516]$} & \shortstack{$0.1098$\\$[-0.1047]$}\\[2pt]
\addlinespace[4pt]
Frank & Gumbel & \shortstack{$0.0540$\\$[+0.0130]$} & \shortstack{$0.0775$\\$[+0.0617]$} & \shortstack{$0.0534$\\$[-0.0217]$} & \shortstack{$0.0576$\\$[+0.0271]$} & \shortstack{$0.1469$\\$[-0.1452]$} & \shortstack{$0.1025$\\$[-0.0965]$}\\[2pt]
Frank & Clayton & \shortstack{$0.0548$\\$[-0.0306]$} & \shortstack{$0.0518$\\$[+0.0083]$} & \shortstack{$0.0739$\\$[-0.0641]$} & \shortstack{$0.0539$\\$[-0.0246]$} & \shortstack{$0.1756$\\$[-0.1751]$} & \shortstack{$0.1346$\\$[-0.1332]$}\\[2pt]
\textbf{Frank} & \textbf{Frank} & \shortstack{$0.0485$\\$[-0.0090]$} & \shortstack{$0.0613$\\$[+0.0376]$} & \shortstack{$0.0589$\\$[-0.0417]$} & \shortstack{$0.0501$\\$[+0.0046]$} & \shortstack{$0.1580$\\$[-0.1570]$} & \shortstack{$0.1154$\\$[-0.1123]$}\\[2pt]
\bottomrule
\end{tabular}
\end{table}
\clearpage

\begin{table}[!htbp]
\centering\small
\caption{Systematic comparison at $n=100$ and $a=q=0.5$. Each entry gives mean absolute error above, with signed bias in square brackets below; both refer to the uncontaminated $\tau_0$. All nine generating/fitted-family combinations and all three contamination fractions are included for each displayed dependence level. Bold family pairs identify correct specification. Every entry uses 200 replications. MCSEs and results for all other evaluated powers and divergence estimators are reported in Tables~\ref{tab:complete-n100-gumbel-gumbel}--\ref{tab:complete-n100-frank-frank}.}
\label{tab:comparison-n100-high}
\setlength{\tabcolsep}{3.5pt}
\renewcommand{\arraystretch}{1.05}
\begin{tabular}{@{}ll*{6}{r}@{}}
\toprule
& & \multicolumn{2}{c}{$\varepsilon=0$} & \multicolumn{2}{c}{$\varepsilon=0.1$} & \multicolumn{2}{c}{$\varepsilon=0.5$}\\
\cmidrule(lr){3-4}\cmidrule(lr){5-6}\cmidrule(lr){7-8}
Generating & Fitted & MCopPDE & EK--$\beta$ & MCopPDE & EK--$\beta$ & MCopPDE & EK--$\beta$\\
\midrule
\multicolumn{8}{@{}l}{\textbf{Kendall tau $\tau_0=0.5$}}\\[3pt]
\textbf{Gumbel} & \textbf{Gumbel} & \shortstack{$0.0439$\\$[-0.0154]$} & \shortstack{$0.0536$\\$[+0.0337]$} & \shortstack{$0.0737$\\$[-0.0680]$} & \shortstack{$0.0488$\\$[-0.0171]$} & \shortstack{$0.2699$\\$[-0.2699]$} & \shortstack{$0.2173$\\$[-0.2173]$}\\[2pt]
Gumbel & Clayton & \shortstack{$0.0653$\\$[-0.0521]$} & \shortstack{$0.0546$\\$[-0.0040]$} & \shortstack{$0.1140$\\$[-0.1118]$} & \shortstack{$0.0770$\\$[-0.0653]$} & \shortstack{$0.3126$\\$[-0.3126]$} & \shortstack{$0.2671$\\$[-0.2671]$}\\[2pt]
Gumbel & Frank & \shortstack{$0.0554$\\$[-0.0450]$} & \shortstack{$0.0455$\\$[+0.0078]$} & \shortstack{$0.0982$\\$[-0.0966]$} & \shortstack{$0.0581$\\$[-0.0443]$} & \shortstack{$0.2888$\\$[-0.2888]$} & \shortstack{$0.2399$\\$[-0.2399]$}\\[2pt]
\addlinespace[4pt]
Clayton & Gumbel & \shortstack{$0.0525$\\$[+0.0348]$} & \shortstack{$0.1130$\\$[+0.1106]$} & \shortstack{$0.0491$\\$[-0.0258]$} & \shortstack{$0.0637$\\$[+0.0476]$} & \shortstack{$0.2474$\\$[-0.2474]$} & \shortstack{$0.1874$\\$[-0.1873]$}\\[2pt]
\textbf{Clayton} & \textbf{Clayton} & \shortstack{$0.0413$\\$[-0.0131]$} & \shortstack{$0.0487$\\$[+0.0315]$} & \shortstack{$0.0740$\\$[-0.0697]$} & \shortstack{$0.0461$\\$[-0.0232]$} & \shortstack{$0.2787$\\$[-0.2787]$} & \shortstack{$0.2296$\\$[-0.2296]$}\\
\bottomrule
\end{tabular}
\end{table}
\clearpage

\section{Complete numerical tables}

\begin{table}[!htbp]
\centering\footnotesize
\caption{Complete results for $n=100$, generating Gumbel and fitted Gumbel. MAE is mean absolute error; bias is signed; SE is the Monte Carlo standard error of MAE in units of $10^{-4}$. Both errors refer to the uncontaminated $\tau_0$. Power denotes $a$ for MCopPDE and $q$ for the EK estimators. Every entry uses 200 replications.}
\label{tab:complete-n100-gumbel-gumbel}
\setlength{\tabcolsep}{3pt}
\renewcommand{\arraystretch}{1.00}
\begin{tabular}{@{}lr*{9}{r}@{}}
\toprule
& & \multicolumn{3}{c}{$\varepsilon=0$} & \multicolumn{3}{c}{$\varepsilon=0.1$} & \multicolumn{3}{c}{$\varepsilon=0.5$}\\
\cmidrule(lr){3-5}\cmidrule(lr){6-8}\cmidrule(lr){9-11}
Estimator & Power & MAE & Bias & SE & MAE & Bias & SE & MAE & Bias & SE\\
\midrule
\multicolumn{11}{@{}l}{\textbf{Kendall tau $\tau_0=0.1$}}\\[1pt]
MCopPDE & 0 & 0.0574 & $-0.0029$ & 29.1 & 0.0589 & $-0.0117$ & 26.0 & 0.0593 & $-0.0383$ & 26.2\\
MCopPDE & 0.1 & 0.0566 & $-0.0016$ & 29.0 & 0.0582 & $-0.0104$ & 25.8 & 0.0588 & $-0.0374$ & 26.1\\
MCopPDE & 0.5 & 0.0546 & $+0.0027$ & 28.7 & 0.0560 & $-0.0060$ & 25.6 & 0.0566 & $-0.0344$ & 25.7\\
MCopPDE & 1 & 0.0537 & $+0.0073$ & 28.0 & 0.0543 & $-0.0019$ & 25.8 & 0.0544 & $-0.0314$ & 25.2\\
EK--$\alpha$ & 0.1 & 0.0802 & $+0.0613$ & 41.7 & 0.0749 & $+0.0514$ & 38.4 & 0.0566 & $+0.0210$ & 31.9\\
\bottomrule
\end{tabular}
\end{table}
\clearpage

\begin{table}[!htbp]
\centering\footnotesize
\caption{Complete results for $n=100$, generating Gumbel and fitted Clayton. MAE is mean absolute error; bias is signed; SE is the Monte Carlo standard error of MAE in units of $10^{-4}$. Both errors refer to the uncontaminated $\tau_0$. Power denotes $a$ for MCopPDE and $q$ for the EK estimators. Every entry uses 200 replications.}
\label{tab:complete-n100-gumbel-clayton}
\setlength{\tabcolsep}{3pt}
\renewcommand{\arraystretch}{1.00}
\begin{tabular}{@{}lr*{9}{r}@{}}
\toprule
& & \multicolumn{3}{c}{$\varepsilon=0$} & \multicolumn{3}{c}{$\varepsilon=0.1$} & \multicolumn{3}{c}{$\varepsilon=0.5$}\\
\cmidrule(lr){3-5}\cmidrule(lr){6-8}\cmidrule(lr){9-11}
Estimator & Power & MAE & Bias & SE & MAE & Bias & SE & MAE & Bias & SE\\
\midrule
\multicolumn{11}{@{}l}{\textbf{Kendall tau $\tau_0=0.1$}}\\[1pt]
MCopPDE & 0 & 0.0550 & $-0.0391$ & 24.1 & 0.0574 & $-0.0451$ & 23.7 & 0.0675 & $-0.0619$ & 21.9\\
MCopPDE & 0.1 & 0.0539 & $-0.0359$ & 24.1 & 0.0565 & $-0.0422$ & 23.6 & 0.0662 & $-0.0598$ & 22.2\\
MCopPDE & 0.5 & 0.0513 & $-0.0234$ & 24.7 & 0.0538 & $-0.0304$ & 23.6 & 0.0618 & $-0.0514$ & 23.4\\
MCopPDE & 1 & 0.0516 & $-0.0079$ & 26.6 & 0.0535 & $-0.0156$ & 24.3 & 0.0576 & $-0.0412$ & 25.0\\
EK--$\alpha$ & 0.1 & 0.0457 & $+0.0104$ & 26.1 & 0.0444 & $+0.0043$ & 22.5 & 0.0396 & $-0.0155$ & 22.4\\
\bottomrule
\end{tabular}
\end{table}
\clearpage

\begin{table}[!htbp]
\centering\footnotesize
\caption{Complete results for $n=100$, generating Gumbel and fitted Frank. MAE is mean absolute error; bias is signed; SE is the Monte Carlo standard error of MAE in units of $10^{-4}$. Both errors refer to the uncontaminated $\tau_0$. Power denotes $a$ for MCopPDE and $q$ for the EK estimators. Every entry uses 200 replications.}
\label{tab:complete-n100-gumbel-frank}
\setlength{\tabcolsep}{3pt}
\renewcommand{\arraystretch}{1.00}
\begin{tabular}{@{}lr*{9}{r}@{}}
\toprule
& & \multicolumn{3}{c}{$\varepsilon=0$} & \multicolumn{3}{c}{$\varepsilon=0.1$} & \multicolumn{3}{c}{$\varepsilon=0.5$}\\
\cmidrule(lr){3-5}\cmidrule(lr){6-8}\cmidrule(lr){9-11}
Estimator & Power & MAE & Bias & SE & MAE & Bias & SE & MAE & Bias & SE\\
\midrule
\multicolumn{11}{@{}l}{\textbf{Kendall tau $\tau_0=0.1$}}\\[1pt]
MCopPDE & 0 & 0.0543 & $-0.0148$ & 26.2 & 0.0562 & $-0.0227$ & 24.5 & 0.0607 & $-0.0458$ & 24.7\\
MCopPDE & 0.1 & 0.0536 & $-0.0134$ & 26.1 & 0.0556 & $-0.0212$ & 24.3 & 0.0601 & $-0.0449$ & 24.8\\
MCopPDE & 0.5 & 0.0518 & $-0.0077$ & 26.4 & 0.0539 & $-0.0155$ & 24.1 & 0.0578 & $-0.0411$ & 25.0\\
MCopPDE & 1 & 0.0520 & $-0.0005$ & 26.7 & 0.0531 & $-0.0088$ & 24.9 & 0.0557 & $-0.0363$ & 25.0\\
EK--$\alpha$ & 0.1 & 0.0657 & $+0.0411$ & 36.0 & 0.0627 & $+0.0325$ & 32.0 & 0.0492 & $+0.0055$ & 28.0\\
\bottomrule
\end{tabular}
\end{table}
\clearpage

\begin{table}[!htbp]
\centering\footnotesize
\caption{Complete results for $n=100$, generating Clayton and fitted Gumbel. MAE is mean absolute error; bias is signed; SE is the Monte Carlo standard error of MAE in units of $10^{-4}$. Both errors refer to the uncontaminated $\tau_0$. Power denotes $a$ for MCopPDE and $q$ for the EK estimators. Every entry uses 200 replications.}
\label{tab:complete-n100-clayton-gumbel}
\setlength{\tabcolsep}{3pt}
\renewcommand{\arraystretch}{1.00}
\begin{tabular}{@{}lr*{9}{r}@{}}
\toprule
& & \multicolumn{3}{c}{$\varepsilon=0$} & \multicolumn{3}{c}{$\varepsilon=0.1$} & \multicolumn{3}{c}{$\varepsilon=0.5$}\\
\cmidrule(lr){3-5}\cmidrule(lr){6-8}\cmidrule(lr){9-11}
Estimator & Power & MAE & Bias & SE & MAE & Bias & SE & MAE & Bias & SE\\
\midrule
\multicolumn{11}{@{}l}{\textbf{Kendall tau $\tau_0=0.1$}}\\[1pt]
MCopPDE & 0 & 0.0634 & $+0.0227$ & 37.1 & 0.0590 & $+0.0110$ & 32.1 & 0.0615 & $-0.0309$ & 25.3\\
MCopPDE & 0.1 & 0.0622 & $+0.0198$ & 35.9 & 0.0585 & $+0.0083$ & 31.0 & 0.0610 & $-0.0321$ & 25.4\\
MCopPDE & 0.5 & 0.0588 & $+0.0098$ & 31.6 & 0.0573 & $-0.0007$ & 27.6 & 0.0605 & $-0.0364$ & 24.8\\
MCopPDE & 1 & 0.0557 & $+0.0006$ & 28.0 & 0.0556 & $-0.0085$ & 25.4 & 0.0608 & $-0.0402$ & 24.4\\
EK--$\alpha$ & 0.1 & 0.1009 & $+0.0903$ & 51.0 & 0.0888 & $+0.0786$ & 47.6 & 0.0665 & $+0.0296$ & 33.1\\
\bottomrule
\end{tabular}
\end{table}
\clearpage

\begin{table}[!htbp]
\centering\footnotesize
\caption{Complete results for $n=100$, generating Clayton and fitted Clayton. MAE is mean absolute error; bias is signed; SE is the Monte Carlo standard error of MAE in units of $10^{-4}$. Both errors refer to the uncontaminated $\tau_0$. Power denotes $a$ for MCopPDE and $q$ for the EK estimators. Every entry uses 200 replications.}
\label{tab:complete-n100-clayton-clayton}
\setlength{\tabcolsep}{3pt}
\renewcommand{\arraystretch}{1.00}
\begin{tabular}{@{}lr*{9}{r}@{}}
\toprule
& & \multicolumn{3}{c}{$\varepsilon=0$} & \multicolumn{3}{c}{$\varepsilon=0.1$} & \multicolumn{3}{c}{$\varepsilon=0.5$}\\
\cmidrule(lr){3-5}\cmidrule(lr){6-8}\cmidrule(lr){9-11}
Estimator & Power & MAE & Bias & SE & MAE & Bias & SE & MAE & Bias & SE\\
\midrule
\multicolumn{11}{@{}l}{\textbf{Kendall tau $\tau_0=0.1$}}\\[1pt]
MCopPDE & 0 & 0.0489 & $-0.0109$ & 25.5 & 0.0473 & $-0.0199$ & 23.8 & 0.0599 & $-0.0515$ & 24.6\\
MCopPDE & 0.1 & 0.0495 & $-0.0097$ & 25.5 & 0.0477 & $-0.0187$ & 23.7 & 0.0596 & $-0.0504$ & 24.8\\
MCopPDE & 0.5 & 0.0515 & $-0.0054$ & 26.1 & 0.0499 & $-0.0146$ & 24.0 & 0.0587 & $-0.0466$ & 25.4\\
MCopPDE & 1 & 0.0541 & $-0.0005$ & 27.7 & 0.0536 & $-0.0099$ & 24.8 & 0.0598 & $-0.0424$ & 24.9\\
EK--$\alpha$ & 0.1 & 0.0552 & $+0.0407$ & 33.6 & 0.0480 & $+0.0321$ & 30.1 & 0.0442 & $-0.0037$ & 21.4\\
\bottomrule
\end{tabular}
\end{table}
\clearpage

\begin{table}[!htbp]
\centering\footnotesize
\caption{Complete results for $n=100$, generating Clayton and fitted Frank. MAE is mean absolute error; bias is signed; SE is the Monte Carlo standard error of MAE in units of $10^{-4}$. Both errors refer to the uncontaminated $\tau_0$. Power denotes $a$ for MCopPDE and $q$ for the EK estimators. Every entry uses 200 replications.}
\label{tab:complete-n100-clayton-frank}
\setlength{\tabcolsep}{3pt}
\renewcommand{\arraystretch}{1.00}
\begin{tabular}{@{}lr*{9}{r}@{}}
\toprule
& & \multicolumn{3}{c}{$\varepsilon=0$} & \multicolumn{3}{c}{$\varepsilon=0.1$} & \multicolumn{3}{c}{$\varepsilon=0.5$}\\
\cmidrule(lr){3-5}\cmidrule(lr){6-8}\cmidrule(lr){9-11}
Estimator & Power & MAE & Bias & SE & MAE & Bias & SE & MAE & Bias & SE\\
\midrule
\multicolumn{11}{@{}l}{\textbf{Kendall tau $\tau_0=0.1$}}\\[1pt]
MCopPDE & 0 & 0.0568 & $+0.0103$ & 31.8 & 0.0533 & $-0.0003$ & 28.1 & 0.0593 & $-0.0380$ & 25.3\\
MCopPDE & 0.1 & 0.0562 & $+0.0084$ & 31.0 & 0.0534 & $-0.0020$ & 27.2 & 0.0593 & $-0.0387$ & 25.3\\
MCopPDE & 0.5 & 0.0548 & $+0.0025$ & 28.7 & 0.0538 & $-0.0072$ & 25.5 & 0.0595 & $-0.0408$ & 25.0\\
MCopPDE & 1 & 0.0542 & $-0.0019$ & 27.1 & 0.0543 & $-0.0107$ & 24.5 & 0.0605 & $-0.0423$ & 24.6\\
EK--$\alpha$ & 0.1 & 0.0825 & $+0.0690$ & 44.9 & 0.0720 & $+0.0583$ & 41.2 & 0.0579 & $+0.0142$ & 27.9\\
\bottomrule
\end{tabular}
\end{table}
\clearpage

\begin{table}[!htbp]
\centering\footnotesize
\caption{Complete results for $n=100$, generating Frank and fitted Gumbel. MAE is mean absolute error; bias is signed; SE is the Monte Carlo standard error of MAE in units of $10^{-4}$. Both errors refer to the uncontaminated $\tau_0$. Power denotes $a$ for MCopPDE and $q$ for the EK estimators. Every entry uses 200 replications.}
\label{tab:complete-n100-frank-gumbel}
\setlength{\tabcolsep}{3pt}
\renewcommand{\arraystretch}{1.00}
\begin{tabular}{@{}lr*{9}{r}@{}}
\toprule
& & \multicolumn{3}{c}{$\varepsilon=0$} & \multicolumn{3}{c}{$\varepsilon=0.1$} & \multicolumn{3}{c}{$\varepsilon=0.5$}\\
\cmidrule(lr){3-5}\cmidrule(lr){6-8}\cmidrule(lr){9-11}
Estimator & Power & MAE & Bias & SE & MAE & Bias & SE & MAE & Bias & SE\\
\midrule
\multicolumn{11}{@{}l}{\textbf{Kendall tau $\tau_0=0.1$}}\\[1pt]
MCopPDE & 0 & 0.0575 & $+0.0018$ & 27.9 & 0.0595 & $-0.0076$ & 26.1 & 0.0666 & $-0.0482$ & 24.0\\
MCopPDE & 0.1 & 0.0569 & $+0.0030$ & 27.8 & 0.0588 & $-0.0066$ & 26.1 & 0.0660 & $-0.0475$ & 24.0\\
MCopPDE & 0.5 & 0.0551 & $+0.0057$ & 27.4 & 0.0569 & $-0.0044$ & 25.6 & 0.0637 & $-0.0448$ & 23.9\\
MCopPDE & 1 & 0.0536 & $+0.0060$ & 27.5 & 0.0549 & $-0.0042$ & 25.9 & 0.0619 & $-0.0429$ & 24.4\\
EK--$\alpha$ & 0.1 & 0.0806 & $+0.0654$ & 40.5 & 0.0761 & $+0.0540$ & 38.7 & 0.0582 & $+0.0048$ & 30.0\\
\bottomrule
\end{tabular}
\end{table}
\clearpage

\begin{table}[!htbp]
\centering\footnotesize
\caption{Complete results for $n=100$, generating Frank and fitted Clayton. MAE is mean absolute error; bias is signed; SE is the Monte Carlo standard error of MAE in units of $10^{-4}$. Both errors refer to the uncontaminated $\tau_0$. Power denotes $a$ for MCopPDE and $q$ for the EK estimators. Every entry uses 200 replications.}
\label{tab:complete-n100-frank-clayton}
\setlength{\tabcolsep}{3pt}
\renewcommand{\arraystretch}{1.00}
\begin{tabular}{@{}lr*{9}{r}@{}}
\toprule
& & \multicolumn{3}{c}{$\varepsilon=0$} & \multicolumn{3}{c}{$\varepsilon=0.1$} & \multicolumn{3}{c}{$\varepsilon=0.5$}\\
\cmidrule(lr){3-5}\cmidrule(lr){6-8}\cmidrule(lr){9-11}
Estimator & Power & MAE & Bias & SE & MAE & Bias & SE & MAE & Bias & SE\\
\midrule
\multicolumn{11}{@{}l}{\textbf{Kendall tau $\tau_0=0.1$}}\\[1pt]
MCopPDE & 0 & 0.0533 & $-0.0353$ & 22.8 & 0.0552 & $-0.0414$ & 24.5 & 0.0714 & $-0.0678$ & 22.7\\
MCopPDE & 0.1 & 0.0521 & $-0.0316$ & 22.8 & 0.0544 & $-0.0380$ & 24.0 & 0.0700 & $-0.0660$ & 23.4\\
MCopPDE & 0.5 & 0.0501 & $-0.0180$ & 23.4 & 0.0527 & $-0.0257$ & 23.5 & 0.0670 & $-0.0589$ & 23.2\\
MCopPDE & 1 & 0.0516 & $-0.0033$ & 25.4 & 0.0537 & $-0.0125$ & 24.6 & 0.0643 & $-0.0500$ & 23.3\\
EK--$\alpha$ & 0.1 & 0.0443 & $+0.0135$ & 23.9 & 0.0455 & $+0.0061$ & 21.9 & 0.0468 & $-0.0266$ & 22.8\\
\bottomrule
\end{tabular}
\end{table}
\clearpage

\begin{table}[!htbp]
\centering\footnotesize
\caption{Complete results for $n=100$, generating Frank and fitted Frank. MAE is mean absolute error; bias is signed; SE is the Monte Carlo standard error of MAE in units of $10^{-4}$. Both errors refer to the uncontaminated $\tau_0$. Power denotes $a$ for MCopPDE and $q$ for the EK estimators. Every entry uses 200 replications.}
\label{tab:complete-n100-frank-frank}
\setlength{\tabcolsep}{3pt}
\renewcommand{\arraystretch}{1.00}
\begin{tabular}{@{}lr*{9}{r}@{}}
\toprule
& & \multicolumn{3}{c}{$\varepsilon=0$} & \multicolumn{3}{c}{$\varepsilon=0.1$} & \multicolumn{3}{c}{$\varepsilon=0.5$}\\
\cmidrule(lr){3-5}\cmidrule(lr){6-8}\cmidrule(lr){9-11}
Estimator & Power & MAE & Bias & SE & MAE & Bias & SE & MAE & Bias & SE\\
\midrule
\multicolumn{11}{@{}l}{\textbf{Kendall tau $\tau_0=0.1$}}\\[1pt]
MCopPDE & 0 & 0.0533 & $-0.0093$ & 25.2 & 0.0560 & $-0.0177$ & 24.2 & 0.0663 & $-0.0540$ & 23.9\\
MCopPDE & 0.1 & 0.0529 & $-0.0078$ & 25.2 & 0.0555 & $-0.0163$ & 24.2 & 0.0661 & $-0.0533$ & 23.7\\
MCopPDE & 0.5 & 0.0519 & $-0.0029$ & 25.2 & 0.0540 & $-0.0121$ & 24.5 & 0.0644 & $-0.0499$ & 23.4\\
MCopPDE & 1 & 0.0519 & $+0.0015$ & 26.0 & 0.0536 & $-0.0081$ & 24.9 & 0.0625 & $-0.0460$ & 23.9\\
EK--$\alpha$ & 0.1 & 0.0667 & $+0.0466$ & 34.7 & 0.0647 & $+0.0364$ & 32.2 & 0.0541 & $-0.0086$ & 26.3\\
\bottomrule
\end{tabular}
\end{table}
\clearpage

\begin{table}[!htbp]
\centering\footnotesize
\caption{Complete results for $n=1000$, generating Gumbel and fitted Gumbel. MAE is mean absolute error; bias is signed; SE is the Monte Carlo standard error of MAE in units of $10^{-4}$. Both errors refer to the uncontaminated $\tau_0$. Power denotes $a$ for MCopPDE and $q$ for the EK estimators. Every entry uses 200 replications.}
\label{tab:complete-n1000-gumbel-gumbel}
\setlength{\tabcolsep}{3pt}
\renewcommand{\arraystretch}{1.00}
\begin{tabular}{@{}lr*{9}{r}@{}}
\toprule
& & \multicolumn{3}{c}{$\varepsilon=0$} & \multicolumn{3}{c}{$\varepsilon=0.1$} & \multicolumn{3}{c}{$\varepsilon=0.5$}\\
\cmidrule(lr){3-5}\cmidrule(lr){6-8}\cmidrule(lr){9-11}
Estimator & Power & MAE & Bias & SE & MAE & Bias & SE & MAE & Bias & SE\\
\midrule
\multicolumn{11}{@{}l}{\textbf{Kendall tau $\tau_0=0.1$}}\\[1pt]
MCopPDE & 0 & 0.0186 & $+0.0019$ & 10.8 & 0.0219 & $-0.0082$ & 11.9 & 0.0518 & $-0.0511$ & 15.7\\
MCopPDE & 0.1 & 0.0182 & $+0.0019$ & 10.6 & 0.0216 & $-0.0083$ & 11.7 & 0.0517 & $-0.0511$ & 15.5\\
MCopPDE & 0.5 & 0.0174 & $+0.0017$ & 9.5 & 0.0206 & $-0.0084$ & 11.2 & 0.0513 & $-0.0509$ & 15.1\\
MCopPDE & 1 & 0.0170 & $+0.0015$ & 8.7 & 0.0196 & $-0.0086$ & 10.8 & 0.0510 & $-0.0508$ & 14.7\\
EK--$\alpha$ & 0.1 & 0.0199 & $+0.0077$ & 11.1 & 0.0210 & $-0.0024$ & 11.5 & 0.0463 & $-0.0453$ & 15.6\\
\bottomrule
\end{tabular}
\end{table}
\clearpage

\begin{table}[!htbp]
\centering\footnotesize
\caption{Complete results for $n=1000$, generating Gumbel and fitted Clayton. MAE is mean absolute error; bias is signed; SE is the Monte Carlo standard error of MAE in units of $10^{-4}$. Both errors refer to the uncontaminated $\tau_0$. Power denotes $a$ for MCopPDE and $q$ for the EK estimators. Every entry uses 200 replications.}
\label{tab:complete-n1000-gumbel-clayton}
\setlength{\tabcolsep}{3pt}
\renewcommand{\arraystretch}{1.00}
\begin{tabular}{@{}lr*{9}{r}@{}}
\toprule
& & \multicolumn{3}{c}{$\varepsilon=0$} & \multicolumn{3}{c}{$\varepsilon=0.1$} & \multicolumn{3}{c}{$\varepsilon=0.5$}\\
\cmidrule(lr){3-5}\cmidrule(lr){6-8}\cmidrule(lr){9-11}
Estimator & Power & MAE & Bias & SE & MAE & Bias & SE & MAE & Bias & SE\\
\midrule
\multicolumn{11}{@{}l}{\textbf{Kendall tau $\tau_0=0.1$}}\\[1pt]
MCopPDE & 0 & 0.0361 & $-0.0355$ & 12.2 & 0.0427 & $-0.0423$ & 13.1 & 0.0704 & $-0.0704$ & 11.7\\
MCopPDE & 0.1 & 0.0335 & $-0.0328$ & 12.1 & 0.0403 & $-0.0398$ & 13.3 & 0.0690 & $-0.0690$ & 12.0\\
MCopPDE & 0.5 & 0.0249 & $-0.0222$ & 11.1 & 0.0317 & $-0.0302$ & 13.3 & 0.0634 & $-0.0634$ & 13.0\\
MCopPDE & 1 & 0.0179 & $-0.0094$ & 9.6 & 0.0240 & $-0.0185$ & 12.0 & 0.0568 & $-0.0566$ & 14.2\\
EK--$\alpha$ & 0.1 & 0.0331 & $-0.0324$ & 11.7 & 0.0393 & $-0.0389$ & 12.7 & 0.0660 & $-0.0660$ & 11.6\\
\bottomrule
\end{tabular}
\end{table}
\clearpage

\begin{table}[!htbp]
\centering\footnotesize
\caption{Complete results for $n=1000$, generating Gumbel and fitted Frank. MAE is mean absolute error; bias is signed; SE is the Monte Carlo standard error of MAE in units of $10^{-4}$. Both errors refer to the uncontaminated $\tau_0$. Power denotes $a$ for MCopPDE and $q$ for the EK estimators. Every entry uses 200 replications.}
\label{tab:complete-n1000-gumbel-frank}
\setlength{\tabcolsep}{3pt}
\renewcommand{\arraystretch}{1.00}
\begin{tabular}{@{}lr*{9}{r}@{}}
\toprule
& & \multicolumn{3}{c}{$\varepsilon=0$} & \multicolumn{3}{c}{$\varepsilon=0.1$} & \multicolumn{3}{c}{$\varepsilon=0.5$}\\
\cmidrule(lr){3-5}\cmidrule(lr){6-8}\cmidrule(lr){9-11}
Estimator & Power & MAE & Bias & SE & MAE & Bias & SE & MAE & Bias & SE\\
\midrule
\multicolumn{11}{@{}l}{\textbf{Kendall tau $\tau_0=0.1$}}\\[1pt]
MCopPDE & 0 & 0.0192 & $-0.0106$ & 10.7 & 0.0255 & $-0.0196$ & 12.6 & 0.0578 & $-0.0575$ & 14.7\\
MCopPDE & 0.1 & 0.0188 & $-0.0101$ & 10.5 & 0.0250 & $-0.0191$ & 12.5 & 0.0574 & $-0.0571$ & 14.6\\
MCopPDE & 0.5 & 0.0175 & $-0.0077$ & 9.6 & 0.0234 & $-0.0169$ & 11.8 & 0.0559 & $-0.0557$ & 14.5\\
MCopPDE & 1 & 0.0167 & $-0.0040$ & 8.8 & 0.0214 & $-0.0136$ & 11.3 & 0.0537 & $-0.0535$ & 14.6\\
EK--$\alpha$ & 0.1 & 0.0177 & $-0.0058$ & 10.4 & 0.0230 & $-0.0147$ & 11.8 & 0.0531 & $-0.0526$ & 14.7\\
\bottomrule
\end{tabular}
\end{table}
\clearpage

\begin{table}[!htbp]
\centering\footnotesize
\caption{Complete results for $n=1000$, generating Clayton and fitted Gumbel. MAE is mean absolute error; bias is signed; SE is the Monte Carlo standard error of MAE in units of $10^{-4}$. Both errors refer to the uncontaminated $\tau_0$. Power denotes $a$ for MCopPDE and $q$ for the EK estimators. Every entry uses 200 replications.}
\label{tab:complete-n1000-clayton-gumbel}
\setlength{\tabcolsep}{3pt}
\renewcommand{\arraystretch}{1.00}
\begin{tabular}{@{}lr*{9}{r}@{}}
\toprule
& & \multicolumn{3}{c}{$\varepsilon=0$} & \multicolumn{3}{c}{$\varepsilon=0.1$} & \multicolumn{3}{c}{$\varepsilon=0.5$}\\
\cmidrule(lr){3-5}\cmidrule(lr){6-8}\cmidrule(lr){9-11}
Estimator & Power & MAE & Bias & SE & MAE & Bias & SE & MAE & Bias & SE\\
\midrule
\multicolumn{11}{@{}l}{\textbf{Kendall tau $\tau_0=0.1$}}\\[1pt]
MCopPDE & 0 & 0.0372 & $+0.0349$ & 14.7 & 0.0257 & $+0.0196$ & 13.0 & 0.0368 & $-0.0353$ & 16.4\\
MCopPDE & 0.1 & 0.0325 & $+0.0292$ & 13.8 & 0.0226 & $+0.0145$ & 12.1 & 0.0391 & $-0.0382$ & 16.5\\
MCopPDE & 0.5 & 0.0204 & $+0.0103$ & 10.6 & 0.0187 & $-0.0025$ & 10.0 & 0.0478 & $-0.0477$ & 16.0\\
MCopPDE & 1 & 0.0189 & $-0.0070$ & 10.1 & 0.0239 & $-0.0179$ & 11.6 & 0.0562 & $-0.0562$ & 15.2\\
EK--$\alpha$ & 0.1 & 0.0438 & $+0.0424$ & 16.0 & 0.0309 & $+0.0268$ & 14.2 & 0.0319 & $-0.0289$ & 15.3\\
\bottomrule
\end{tabular}
\end{table}
\clearpage

\begin{table}[!htbp]
\centering\footnotesize
\caption{Complete results for $n=1000$, generating Clayton and fitted Clayton. MAE is mean absolute error; bias is signed; SE is the Monte Carlo standard error of MAE in units of $10^{-4}$. Both errors refer to the uncontaminated $\tau_0$. Power denotes $a$ for MCopPDE and $q$ for the EK estimators. Every entry uses 200 replications.}
\label{tab:complete-n1000-clayton-clayton}
\setlength{\tabcolsep}{3pt}
\renewcommand{\arraystretch}{1.00}
\begin{tabular}{@{}lr*{9}{r}@{}}
\toprule
& & \multicolumn{3}{c}{$\varepsilon=0$} & \multicolumn{3}{c}{$\varepsilon=0.1$} & \multicolumn{3}{c}{$\varepsilon=0.5$}\\
\cmidrule(lr){3-5}\cmidrule(lr){6-8}\cmidrule(lr){9-11}
Estimator & Power & MAE & Bias & SE & MAE & Bias & SE & MAE & Bias & SE\\
\midrule
\multicolumn{11}{@{}l}{\textbf{Kendall tau $\tau_0=0.1$}}\\[1pt]
MCopPDE & 0 & 0.0146 & $+0.0007$ & 7.4 & 0.0166 & $-0.0109$ & 9.2 & 0.0520 & $-0.0520$ & 12.9\\
MCopPDE & 0.1 & 0.0147 & $+0.0007$ & 7.5 & 0.0167 & $-0.0109$ & 9.2 & 0.0520 & $-0.0520$ & 13.0\\
MCopPDE & 0.5 & 0.0156 & $+0.0004$ & 7.8 & 0.0176 & $-0.0111$ & 9.5 & 0.0521 & $-0.0521$ & 13.7\\
MCopPDE & 1 & 0.0168 & $-0.0000$ & 8.7 & 0.0192 & $-0.0115$ & 10.1 & 0.0523 & $-0.0523$ & 14.7\\
EK--$\alpha$ & 0.1 & 0.0152 & $+0.0051$ & 7.4 & 0.0147 & $-0.0064$ & 8.4 & 0.0472 & $-0.0472$ & 12.7\\
\bottomrule
\end{tabular}
\end{table}
\clearpage

\begin{table}[!htbp]
\centering\footnotesize
\caption{Complete results for $n=1000$, generating Clayton and fitted Frank. MAE is mean absolute error; bias is signed; SE is the Monte Carlo standard error of MAE in units of $10^{-4}$. Both errors refer to the uncontaminated $\tau_0$. Power denotes $a$ for MCopPDE and $q$ for the EK estimators. Every entry uses 200 replications.}
\label{tab:complete-n1000-clayton-frank}
\setlength{\tabcolsep}{3pt}
\renewcommand{\arraystretch}{1.00}
\begin{tabular}{@{}lr*{9}{r}@{}}
\toprule
& & \multicolumn{3}{c}{$\varepsilon=0$} & \multicolumn{3}{c}{$\varepsilon=0.1$} & \multicolumn{3}{c}{$\varepsilon=0.5$}\\
\cmidrule(lr){3-5}\cmidrule(lr){6-8}\cmidrule(lr){9-11}
Estimator & Power & MAE & Bias & SE & MAE & Bias & SE & MAE & Bias & SE\\
\midrule
\multicolumn{11}{@{}l}{\textbf{Kendall tau $\tau_0=0.1$}}\\[1pt]
MCopPDE & 0 & 0.0262 & $+0.0216$ & 11.5 & 0.0187 & $+0.0079$ & 10.2 & 0.0419 & $-0.0415$ & 15.7\\
MCopPDE & 0.1 & 0.0235 & $+0.0177$ & 10.7 & 0.0177 & $+0.0043$ & 9.6 & 0.0438 & $-0.0435$ & 15.6\\
MCopPDE & 0.5 & 0.0176 & $+0.0052$ & 9.1 & 0.0179 & $-0.0068$ & 9.6 & 0.0498 & $-0.0497$ & 15.2\\
MCopPDE & 1 & 0.0177 & $-0.0056$ & 9.5 & 0.0223 & $-0.0165$ & 11.1 & 0.0550 & $-0.0550$ & 14.9\\
EK--$\alpha$ & 0.1 & 0.0311 & $+0.0282$ & 13.0 & 0.0214 & $+0.0142$ & 11.3 & 0.0369 & $-0.0361$ & 15.4\\
\bottomrule
\end{tabular}
\end{table}
\clearpage

\begin{table}[!htbp]
\centering\footnotesize
\caption{Complete results for $n=1000$, generating Frank and fitted Gumbel. MAE is mean absolute error; bias is signed; SE is the Monte Carlo standard error of MAE in units of $10^{-4}$. Both errors refer to the uncontaminated $\tau_0$. Power denotes $a$ for MCopPDE and $q$ for the EK estimators. Every entry uses 200 replications.}
\label{tab:complete-n1000-frank-gumbel}
\setlength{\tabcolsep}{3pt}
\renewcommand{\arraystretch}{1.00}
\begin{tabular}{@{}lr*{9}{r}@{}}
\toprule
& & \multicolumn{3}{c}{$\varepsilon=0$} & \multicolumn{3}{c}{$\varepsilon=0.1$} & \multicolumn{3}{c}{$\varepsilon=0.5$}\\
\cmidrule(lr){3-5}\cmidrule(lr){6-8}\cmidrule(lr){9-11}
Estimator & Power & MAE & Bias & SE & MAE & Bias & SE & MAE & Bias & SE\\
\midrule
\multicolumn{11}{@{}l}{\textbf{Kendall tau $\tau_0=0.1$}}\\[1pt]
MCopPDE & 0 & 0.0193 & $+0.0105$ & 10.6 & 0.0181 & $-0.0012$ & 10.1 & 0.0451 & $-0.0444$ & 16.4\\
MCopPDE & 0.1 & 0.0189 & $+0.0097$ & 10.2 & 0.0179 & $-0.0019$ & 9.9 & 0.0455 & $-0.0450$ & 16.3\\
MCopPDE & 0.5 & 0.0174 & $+0.0061$ & 9.0 & 0.0176 & $-0.0054$ & 9.4 & 0.0477 & $-0.0474$ & 15.6\\
MCopPDE & 1 & 0.0165 & $+0.0009$ & 8.1 & 0.0184 & $-0.0103$ & 9.7 & 0.0505 & $-0.0504$ & 15.1\\
EK--$\alpha$ & 0.1 & 0.0221 & $+0.0163$ & 11.7 & 0.0185 & $+0.0047$ & 10.5 & 0.0395 & $-0.0384$ & 16.1\\
\bottomrule
\end{tabular}
\end{table}
\clearpage

\begin{table}[!htbp]
\centering\footnotesize
\caption{Complete results for $n=1000$, generating Frank and fitted Clayton. MAE is mean absolute error; bias is signed; SE is the Monte Carlo standard error of MAE in units of $10^{-4}$. Both errors refer to the uncontaminated $\tau_0$. Power denotes $a$ for MCopPDE and $q$ for the EK estimators. Every entry uses 200 replications.}
\label{tab:complete-n1000-frank-clayton}
\setlength{\tabcolsep}{3pt}
\renewcommand{\arraystretch}{1.00}
\begin{tabular}{@{}lr*{9}{r}@{}}
\toprule
& & \multicolumn{3}{c}{$\varepsilon=0$} & \multicolumn{3}{c}{$\varepsilon=0.1$} & \multicolumn{3}{c}{$\varepsilon=0.5$}\\
\cmidrule(lr){3-5}\cmidrule(lr){6-8}\cmidrule(lr){9-11}
Estimator & Power & MAE & Bias & SE & MAE & Bias & SE & MAE & Bias & SE\\
\midrule
\multicolumn{11}{@{}l}{\textbf{Kendall tau $\tau_0=0.1$}}\\[1pt]
MCopPDE & 0 & 0.0282 & $-0.0274$ & 10.9 & 0.0359 & $-0.0354$ & 11.5 & 0.0643 & $-0.0643$ & 12.3\\
MCopPDE & 0.1 & 0.0258 & $-0.0246$ & 10.6 & 0.0336 & $-0.0329$ & 11.4 & 0.0629 & $-0.0629$ & 12.5\\
MCopPDE & 0.5 & 0.0189 & $-0.0143$ & 9.1 & 0.0257 & $-0.0237$ & 11.0 & 0.0578 & $-0.0578$ & 13.4\\
MCopPDE & 1 & 0.0158 & $-0.0034$ & 8.0 & 0.0195 & $-0.0140$ & 10.1 & 0.0524 & $-0.0523$ & 14.6\\
EK--$\alpha$ & 0.1 & 0.0261 & $-0.0250$ & 10.3 & 0.0332 & $-0.0325$ & 11.0 & 0.0600 & $-0.0600$ & 12.0\\
\bottomrule
\end{tabular}
\end{table}
\clearpage

\begin{table}[!htbp]
\centering\footnotesize
\caption{Complete results for $n=1000$, generating Frank and fitted Frank. MAE is mean absolute error; bias is signed; SE is the Monte Carlo standard error of MAE in units of $10^{-4}$. Both errors refer to the uncontaminated $\tau_0$. Power denotes $a$ for MCopPDE and $q$ for the EK estimators. Every entry uses 200 replications.}
\label{tab:complete-n1000-frank-frank}
\setlength{\tabcolsep}{3pt}
\renewcommand{\arraystretch}{1.00}
\begin{tabular}{@{}lr*{9}{r}@{}}
\toprule
& & \multicolumn{3}{c}{$\varepsilon=0$} & \multicolumn{3}{c}{$\varepsilon=0.1$} & \multicolumn{3}{c}{$\varepsilon=0.5$}\\
\cmidrule(lr){3-5}\cmidrule(lr){6-8}\cmidrule(lr){9-11}
Estimator & Power & MAE & Bias & SE & MAE & Bias & SE & MAE & Bias & SE\\
\midrule
\multicolumn{11}{@{}l}{\textbf{Kendall tau $\tau_0=0.1$}}\\[1pt]
MCopPDE & 0 & 0.0161 & $-0.0014$ & 8.6 & 0.0190 & $-0.0119$ & 10.5 & 0.0508 & $-0.0506$ & 15.3\\
MCopPDE & 0.1 & 0.0159 & $-0.0013$ & 8.5 & 0.0189 & $-0.0119$ & 10.3 & 0.0509 & $-0.0507$ & 15.2\\
MCopPDE & 0.5 & 0.0158 & $-0.0012$ & 8.0 & 0.0187 & $-0.0118$ & 9.9 & 0.0511 & $-0.0510$ & 15.0\\
MCopPDE & 1 & 0.0159 & $-0.0011$ & 7.9 & 0.0187 & $-0.0119$ & 9.7 & 0.0511 & $-0.0511$ & 15.0\\
EK--$\alpha$ & 0.1 & 0.0163 & $+0.0037$ & 9.1 & 0.0174 & $-0.0068$ & 9.9 & 0.0459 & $-0.0455$ & 15.2\\
\bottomrule
\end{tabular}
\end{table}
\clearpage